\documentclass[11pt]{article}
\usepackage{amsmath,amssymb}
\usepackage{amsthm}
\usepackage{mathtools}
\usepackage[noend]{algorithmic}
\usepackage[ruled,vlined]{algorithm2e}
\usepackage{url}
\usepackage{fullpage}
\usepackage{makeidx}
\usepackage{enumerate}
\usepackage[top=1in, bottom=1.25in, left=1in, right=1in]{geometry}
\usepackage{graphicx,float,psfrag,epsfig,caption}
\usepackage[usenames,dvipsnames,svgnames,table]{xcolor}
\definecolor{darkgreen}{rgb}{0.0,0,0.9}
\usepackage[pagebackref,letterpaper=true,colorlinks=true,pdfpagemode=none,citecolor=OliveGreen,linkcolor=BrickRed,urlcolor=BrickRed,pdfstartview=FitH]{hyperref}
\usepackage{epstopdf}
\usepackage{color}
\usepackage{xr}
\usepackage{subfig}
\usepackage{caption}
\usepackage{graphicx}
\usepackage[utf8]{inputenc}

\usepackage{scalerel,stackengine}

\stackMath
\newcommand\reallywidehat[1]{%
\savestack{\tmpbox}{\stretchto{%
  \scaleto{%
    \scalerel*[\widthof{\ensuremath{#1}}]{\kern.1pt\mathchar"0362\kern.1pt}%
    {\rule{0ex}{\textheight}}
  }{\textheight}%
}{2.4ex}}%
\stackon[-6.9pt]{#1}{\tmpbox}%
}
\usepackage[mathscr]{euscript}

\DeclareSymbolFont{rsfs}{U}{rsfs}{m}{n}
\DeclareSymbolFontAlphabet{\mathscrsfs}{rsfs}

\numberwithin{equation}{section}

\newtheoremstyle{myexample} 
    {\topsep}                    
    {\topsep}                    
    {\rm }                   
    {}                           
    {\bf }                   
    {.}                          
    {.5em}                       
    {}  

\newtheoremstyle{myremark} 
    {\topsep}                    
    {\topsep}                    
    {\rm}                        
    {}                           
    {\bf}                        
    {.}                          
    {.5em}                       
    {}  

\newtheorem{claim}{Claim}[section]

\newtheorem{theorem}{Theorem}
\newtheorem{proposition}[claim]{Proposition}

\newtheorem{definition}[claim]{Definition}

\theoremstyle{myremark}
\newtheorem{remark}{Remark}[section]

\theoremstyle{myremark}

\theoremstyle{myexample}

\definecolor{darkgreen}{rgb}{0.0, 0.5, 0.0}

\newcommand{\bea}{\begin{eqnarray}}
\newcommand{\eea}{\end{eqnarray}}
\newcommand{\<}{\langle}
\renewcommand{\>}{\rangle}
\newcommand{\E}{{\mathbb E}}

\definecolor{darkblue}{rgb}{0, 0, 0.5}

\def\PL{{\rm PL}_2}
\def\sF{{\sf F}}
\def\sV{{\sf V}}
\def\free{\phi}

\def\TT{{\mathbb{T}}}

\def\ovp{\overline{p}}
\def\SBM{{\sf SBM}}

\def\Info{{\rm I}}
\def\snew{\mbox{\tiny \rm new}}

\def\sTV{\mbox{\tiny \rm TV}}

\def\salg{\mbox{\tiny \rm alg}}

\def\sSDP{\mbox{\tiny \rm SDP}}
\def\sSK{\mbox{\tiny \rm SK}}

\def\Unif{{\sf Unif}}

\def\eps{{\varepsilon}}

\def\bPi{{\boldsymbol{\Pi}}}
\def\id{{\boldsymbol{I}}}

\def\cuC{\mathscrsfs{C}}
\def\cuD{\mathscrsfs{D}}

\def\cuF{\mathscrsfs{F}}

\def\cuE{\mathscrsfs{E}}
\def\cuU{\mathscrsfs{U}}
\def\cuL{\mathscrsfs{L}}
\def\cuS{\mathscrsfs{S}}
\def\cuZ{\mathscrsfs{Z}}

\def\osigma{\overline{\sigma}}

\def\OPT{{\sf OPT}}
\def\ALG{{\sf ALG}}
\def\SDP{{\sf SDP}}
\def\UB{{\sf UB}}

\def\sBayes{\mbox{\tiny\rm Bayes}}

\def\htheta{\hat{\theta}}

\def\Yeff{Y_{\mbox{\tiny\rm eff}}} 

\def\bm{{\boldsymbol{m}}}
\def\bTheta{{\boldsymbol{\Theta}}}
\def\btheta{{\boldsymbol{\theta}}}
\def\hbtheta{\hat{\boldsymbol{\theta}}}

\def\bsigma{{\boldsymbol{\sigma}}}

\def\hsb{\hat{\sf b}}

\newcommand{\indep}{\perp \!\!\! \perp}

\def\sP{{\sf P}}

\def\bC{{\boldsymbol{C}}}
\def\bQ{{\boldsymbol{Q}}}

\def\omu{{\overline{\mu}}}

\def\GOE{{\sf GOE}}

\def\hbQ{{\boldsymbol{\widehat{Q}}}}

\def\bg{{\boldsymbol{g}}}

\def\bzero{{\mathbf 0}}

\def\cC{{\mathcal C}}
\def\cX{{\mathcal X}}

\def\op{\mbox{\tiny\rm op}}

\def\naturals{{\mathbb N}}
\def\reals{{\mathbb R}}

\def\normal{{\sf N}}

\def\sT{{\sf T}}

\def\bz{{\boldsymbol{z}}}
\def\bx{{\boldsymbol{x}}}

\def\bA{\boldsymbol{A}}

\def\bX{\boldsymbol{X}}

\def\de{{\rm d}}

\def\bX{\boldsymbol{X}}
\def\bY{\boldsymbol{Y}}

\def\prob{{\mathbb P}}

\def\E{{\mathbb E}}

\def\<{\langle}
\def\>{\rangle}
\def\Tr{{\sf Tr}}

\def\Ball{{\sf B}}
\def\argmax{{\arg\!\max}}
\def\sign{{\rm sign}}

\def\conv{{\rm conv}}

\def\cL{{\cal L}}

\def\S{{\mathbb S}}

\def\by{{\boldsymbol{y}}}
\def\bw{{\boldsymbol{w}}}

\def\bnu{{\boldsymbol \nu}}

\def\b0{{\boldsymbol{0}}}

\def\bfone{{\boldsymbol 1}}

\def\bG{{\boldsymbol G}}

\def\bm{{\boldsymbol m}}

\def\bs{{\boldsymbol s}}

\DeclareMathOperator*{\plim}{p-lim}

\def\cA{{\mathcal A}}

\def\cS{{\mathcal S}}

\title{Spin Glass Concepts\\ in Computer Science, Statistics, and Learning}

\author{Andrea Montanari\thanks{Department of Statistics and Department of Mathematics,
Stanford University}}

\begin{document}

\maketitle

\begin{abstract}
Spin glass theory studies the structure of sublevel sets and minima (or near-minima) of 
certain classes of random functions in high dimension. 
Near-minima of random functions also play an important role in 
high-dimensional statistics and statistical learning, where minimizing 
the empirical risk (which is a random function of the model parameters) 
is the method of choice for learning a statistical model from noisy data.
Finally, near-minima of random functions are obviously central to average-case
analysis of optimization algorithms.
Computer science, statistics, and machine learning naturally 
lead to questions that are traditionally not addressed within 
physics and mathematical physics. I will try to explain how ideas from 
spin glass theory have seeded recent developments in these fields. 

(This article was written on the occasion of the 2024 Abel Prize to Michel Talagrand.)
\end{abstract}

\tableofcontents

\section{Introduction}

The distribution of the maximum of a collection of (possibly dependent) random 
variables is one of the most classical questions in probability theory 
\cite{FisherTippett1928,Frechet1927,Gnedenko1943}. Upper bounds on the maximum of 
certain stochastic processes were among Kolmogorv's key achievements, 
and instrumental in proving H\"older regularity of stochastic processes
\cite{Kolmogorov1933,Kolmogorov1950}.
Regularity of stochastic processes was one of Michel Talagrand's (MT)
early interests, eventually leading to his powerful generic chaining theory 
\cite{talagrand1987regularity,talagrand2014upper}.

Maxima (or --equivalently-- minima) of random functions
also  play a prominent role in optimization theory (where the random function models an 
objective to be optimized), in statistical physics (where the random function models 
the Hamiltonian of a disordered physical system), and in machine learning (where the 
random function models the deviation between an empirical risk function and its expectation).

MT was a pioneer in recognizing the deep mathematical unity of these domains,
and developing the probabilistic foundations for them 
\cite{rhee1987martingale,rhee1988exact,talagrand1998sherrington,TalagrandVolI,TalagrandVolII}.
A distinct subset of his work stemmed from the desire to understand and make rigorous
the physicists' theory of mean-field spin glasses, as summarized by M\'ezard, Parisi, Virasoro in
\cite{SpinGlass}. 
A large group of brilliant mathematicians contributed to this effort, and 
we refer to \cite{TalagrandVolI,TalagrandVolII,bovier2006statistical,panchenko2013sherrington} for
surveys.

In addition to confirming and elucidating the content of statistical physics ideas, 
this work provided the foundation for extensions of mathematical spin glass theory
beyond  physics. The present article provides a short summary of some of these developments
and their motivations.

As a running example, I will consider a random optimization problem that 
plays a central role in spin glass theory and in MT's work:
\begin{align}
\mbox{maximize}&\;\;\;\; H(\bsigma) =\frac{1}{2} \<\bsigma,\bA\bsigma\> \,,\label{eq:Hsigma}\\
\mbox{subject to}&\;\;\; \;\bsigma \in \Sigma_n\, .
\end{align}
Here $\Sigma_n\subseteq \reals^n$ is a constraint set (the space of `spin configurations' 
in physics parlance.)
A prototypical example is the Ising/hypercube case $\Sigma_n=\{+1,-1\}^n$.
Further, $\bA\in\reals^{n\times n}$ is a symmetric random matrix  and we will consider a few different 
models for its distribution.

More generally, mean field spin glasses are random functions $H:\Sigma_n\to \reals$, where 
$\Sigma_n\subseteq \reals^n$
is a high dimensional product space\footnote{Even more generally,
one can consider $\Sigma_n$ to be a subset of a product space which is symmetric under 
permutations of coordinates, i.e. $\bx\in \Sigma_n\Rightarrow \bPi\bx \in\Sigma_n$ 
for any permutation matrix $\bPi$.}, and $H$ has the property that, for every $k$,
and any $\bsigma_1,\dots,\bsigma_k\in\Sigma_n$ the joint distribution of
$H(\bsigma_1),\dots ,H(\bsigma_k)$  depends on  $\bsigma_1,\dots,\bsigma_k$ only via their
joint empirical distribution\footnote{Also this condition can be somewhat generalized to make the distribution
to depend also on some fixed additional vector $\bw\in\reals^n$, via 
 $n^{-1}\sum_{i\le n}\delta_{\sigma_{1,i}\dots \sigma_{k,i},w_i}$.}
 $n^{-1}\sum_{i\le n}\delta_{\sigma_{1,i}\dots \sigma_{k,i}}$.

In the rest of this introduction, I will describe a few of the
questions that arise from computer science and high-dimensional statistics.
 Section \ref{sec:Parisi} reviews Parisi's formula, whose proof was
 one of MT's main achievements.
 Section \ref{sec:CS-Stats} outlines some of the approaches developed 
in computer science and statistics to address the questions raised 
in the introduction. Section \ref{sec:MP} discusses message passing algorithms,
which play an important role at the nexus between computer science, 
machine learning and physics, and Section \ref{sec:AMP} focuses on a 
version of these algorithms that are best suited for mean field models on 
dense graphs. In Section \ref{sec:Algo} we close the loop by showing that 
the analysis of approximate message passing algorithms yields an 
algorithmic version of Parisi's formula.
The ideas surveyed in this paper have been applied to
a significanlty broader set of problems than what we can 
discuss in this paper. Section \ref{sec:Broader} provides a few pointers
to the literature on related problems.

\subsection{Motivations and questions}
\label{sec:Motivations}

I will recall the motivation for studying
the optimization problem \eqref{eq:Hsigma} in various domains.
 While  \eqref{eq:Hsigma} is an historically important example,
the same questions are interesting for a much broader set of problems, 
as discussed in Section \ref{sec:Broader}.

\paragraph{Statistical physics.} The function $H(\bsigma)$ in Eq.~\eqref{eq:Hsigma}
can be regarded as the Hamiltonian of a spin model, and in particular a Ising model when 
$\Sigma_n=\{+1,-1\}^n$.

Spin models are used to model the interactions of 
molecular magnetic moments in a magnetic material \cite{blundell2001magnetism}.
Each Ising spin
$\sigma_i\in\{+1,-1\}$ corresponds to the alignment of a molecular moment with a special
axis. Interactions between magnetic moments are of quantum-mechanical nature, 
but often this can be abstracted away and the 
 properties of such a system at equilibrium are determined by 
the (classical) Gibbs measure 
\begin{align}
\mu(\de\bsigma) = \frac{1}{Z_{n}(\beta,h)} e^{\beta H(\bsigma)+ h \<\bfone,\bsigma\>}\, \nu_0(\de\bsigma)\, . \label{eq:FirstGibbs}
\end{align}
where $\beta\ge 0$ is the inverse temperature, $h\in\reals^n$ is a `magnetic field,' and $\nu_0$ is a reference measure on $\Sigma_n$ (typically, the uniform measure over
$\Sigma_n$). We will always normalize $\bA$ in  Eq.~\eqref{eq:Hsigma} so that $\|\bA\|_{\op} = \Theta(1)$
and hence $\max_{\bsigma} H(\bsigma) = O(n)$. 
(The sign convention in physics is to take the Gibbs measure to be $\mu(\de\bsigma) \propto
\exp\{-\beta H(\bsigma)\}$: obviously we can flip the sign of $H$ and hence
exchange the role of maxima and minima.)

For large\footnote{The discussion here is heuristic, but this claim holds under many interesting choices of
 $\bA$, for most $\beta$.} $n$, the measure $\mu$ is roughly equivalent to the uniform measure on 
a superlevel set $\{\bsigma\in\Sigma_n:\, H(\bsigma)\ge n\, e(\beta)\}$ for a certain deterministic
function  $\beta\mapsto e(\beta)$. Hence, the geometry of such superlevel sets 
plays a crucial role in the properties of the Gibbs measure $\mu$. 
For large $\beta$, thermodynamic properties are controlled by the maxima and near maxima of $H(\,\cdot\,)$.

In the classical Ising model, $\bA$ is the adjecency matrix of a regular $d$-dimensional grid. 
Since the seventies, physicists became interested in disordered magnetic materials,
in which magnetic moments are localized at random position in an alloy, rather than on a regular lattice. 
These systems were  modeled by Hamiltonians of the form \eqref{eq:Hsigma} where the couplings 
$(A_{ij})_{i,j\le n}$ are taken to be random variables \cite{edwards1975theory}. 
The Gibbs measure $\mu(\,\cdot\,) = \mu_{\bA}(\,\cdot\,)$ 
is thus a random measure (measurable on $\bA$) and we are interested in its typical properties.

Sherrington and Kirkpatrick (SK) \cite{sherrington1975solvable} took the step of defining what is 
arguably the simplest spin
glass model of this class. Namely, they studied the case in which $\bA$ is a matrix from the Gaussian Orthogonal 
Ensemble (we will write $\bA\sim\GOE(n)$). Namely $(A_{ij})_{i\le j\le n}$
are independent  random variables with $A_{ii}\sim\normal(0,2/n)$ and  $A_{ij}\sim\normal(0,1/n)$
for $i<j$. This is the celebrated SK model.

\paragraph{Computer science.} In the case $\Sigma_n=\{+1,-1\}^n$,
\eqref{eq:Hsigma} is a very classical optimization problem known as `binary quadratic optimization.'  
I will denote by $\OPT_n=\OPT_n(\bA)$ its value (rescaled by $n$ for convenience)
\begin{align}
\OPT_n(\bA):=\max_{\bsigma\in\Sigma_n}\, \frac{1}{2n} \<\bsigma,\bA\bsigma\>\, .\label{eq:OPT-Def}
\end{align}
In physics language, this is the `ground state energy' of the Hamiltonian $H(\bsigma)$ 
(up to a sign flip).

Unlike in physics, the questions of interest in computer science are of algorithmic nature.
The most important such question is whether there  exists of a polynomial time algorithm
$\bA\mapsto \bsigma^{\salg}(\bA)\in\Sigma_{n}$ to compute an optimizer or near 
optimizer of this problem.
Namely we would like such an algorithm to achieve
\begin{align}
\frac{1}{2n} \<\bsigma^{\salg},\bA\bsigma^{\salg}\>\ge(1-\eps)\cdot\OPT_n(\bA)\, ,\label{eq:ApproxOPT}
\end{align}
for $\eps$ as small a constant as possible (independent of $n$). An algorithm achieving
this goal is referred to as a $(1-\eps)$-approximation algorithm. Note that, for any fixed $n$,
the exact optimizer can be found by exhaustive enumeration\footnote{Throughout this paper,
we deliberately ignore  the fact that the matrix $A_{ij}$ must be rounded to a finite number of digits
in a finite model of computation. These subtleties can be easily addressed in the present context.}
in time $O(2^n)$.
Hence, in order establish the existence of a $(1-\eps)$-polynomial time approximation algorithm, 
it is only necessary to consider $n\ge n_0$
for a large enough $n_0$. In other words, the question of approximation is really 
a question about large instances. As discussed in Section \ref{sec:Algo}, we do not expect such an 
approximation algorithm to exist, in general (except for very 
poor approximation ratio, i.e.
 $1-\eps$ vanishing with $n$).

Specific subclasses of problem \eqref{eq:Hsigma} played a prominent role in the history of
combinatorial optimization. If $\bA = \bA_G$ is the adjacency matrix\footnote{Namely
$A_{G,ij}=1$ if $(i,j)\in E$ and $A_{G,ij}=0$ otherwise.} of a graph $G=(V=[n],E)$
and $\Sigma_n = \{+1,-1\}^n\cap\{\bx:\, \<\bx,\bfone\>=0\}$, then \eqref{eq:Hsigma}
is the minimum bisection problem (`partition the vertices in two equal sets as to minimize
the number of edges across'). If $\bA = -\bA_G$ and $\Sigma_n = \{+1,-1\}^n$, then 
\eqref{eq:Hsigma} is the max-cut problem (`partition the vertices in two sets as to maximize
the number of edges across'). Both of these sub-problems are NP-complete (see 
Section \ref{sec:Algo}) and in fact among the first ones to motivate the study of
probabilistic models for $\bA$ (i.e., for the graph $G$) \cite{bui1987graph,condon2001algorithms}, as well
as the first breakthroughs in approximation algorithms \cite{goemans1995improved}.

\paragraph{High-dimensional statistics.} A number of modern data-sets can be naturally 
recast as matrices (or higher-order tensors). For instance, the structure of a undirected network.
i.e. a simple graph $G=(V,E)$ with  vertex set $V$ and  edge set $E$, 
can be encoded by its adjacency matrix $\bA_G$. 
An important problem is the one of identifying latent low-dimensional 
structures in such matrix or graph data. The 
balanced two-communities stochastic block
model \cite{holland1983stochastic,abbe2018community} $\SBM(n,p,q)$ is the simplest 
statistical model for this problem. 
To generate $G=(V,E)\sim \SBM(n,p,q)$, we first draw $\btheta\sim\Unif(\Sigma_n)$, 
$\Sigma_n:=\{\bx\in \{+1,-1\}^n:\; \<\bx,\bfone\>=0\}$, encoding the community structure
(for simplicity, we assume $n$ even), 
and then sample edges $(i,j)$ conditionally independently given $\btheta$:
\begin{align}\label{eq:SBM-def}
\prob\big(A_{G,ij}=1\big|\btheta\big)=
\begin{cases}
p & \mbox{ if $\theta_i\theta_j=+1$},\\
q & \mbox{ if $\theta_i\theta_j=-1$},
\end{cases}
\end{align}
It is understood that $\bA_G = \bA_G^{\sT}$ with $A_{G,ij}\in\{0,1\}$.
Given a single realization of the matrix $\bA_G$, we would like to estimate the community structure
$\btheta$ (which is of course identifiable only up to a global flip $\btheta\mapsto -\btheta$).

For $p>q$ edges are more likely to connect vertices $i,j$ in the same 
community  (i.e., with $\theta_i=\theta_j$) rather than vertices in different communities
(with $\theta_i\neq\theta_j$). It is then easy to check that the maximum likelihood 
estimator is the solution of problem \eqref{eq:Hsigma}, namely
\begin{align}\label{eq:MLE}
\hbtheta(G) := \argmax\big\{\<\bsigma,\bA_G\bsigma\>\; :\;\;\bsigma\in \Sigma_n\, \big\}\, .
\end{align}
This is just a minimum bisection of graph $G$. 

The value of this optimization problem (the `ground state energy')
$\OPT_n(\bA_G;\Sigma_n)$ has a statistical application.
Namely, it can be used to `detect' the presence of a latent 
community structure. 
For instance, assume we are given a graph $G\sim \SBM(n,\ovp+\Delta,\ovp-\Delta)$ and we know $\ovp$
(this is reasonable because $\binom{n}{2}\ovp$ is the expected number of edges, and the actual number of edges
concentrates around its expectation). We are asked to decide whether $\Delta=0$ or $\Delta\ge 
\underline{\Delta}_{n}$. The corresponding mathematical problem is: 
what is the smallest sequence of constants
 $\underline{\Delta}_{n}$ such that 
 distinguishing between the two hypotheses $\Delta=0$ and $\Delta\ge 
\underline{\Delta}_{n}$
 is possible with error  probability bounded by some $p_0$?
 One interesting approach is to compare the ground state energy 
 $\OPT_n(\bA_G;\Sigma_n)$ to its expectation 
 $\E_0\OPT_n(\bA_G;\Sigma_n)$
 under the `null model'
 $\bA_G\sim\SBM(n,\ovp,\ovp)$.

An even more interesting question is the one of estimating $\btheta$, and
in particular understanding the properties of the MLE \eqref{eq:MLE}:
$(i)$~What is the statistical accuracy of $\hbtheta(G)$? $(ii)$~Is this the optimal accuracy 
achieved by any estimator?
$(iii)$~Can we compute $\hbtheta(G)$ (or a superior estimator) in polynomial time?
Accuracy can be measured via the normalized inner product (overlap)
\begin{align}
Q_{\hbtheta} = \frac{1}{n}\E\big\{|\<\hbtheta(G),\btheta\>|\big\}\, ,
\end{align}
which is the absolute value of the 
fraction of vertices whose community membership is estimated correctly,
minus the fraction that is estimated incorrectly.

\paragraph{Sampling.} The focus of this short review will be on
the application of spin glass techniques to optimization 
(e.g., problem \eqref{eq:Hsigma}) and statistical estimation
(e.g., estimating $\btheta$ from data \eqref{eq:SBM-def}).
An equally important problem in computer science, statistics, and
 machine learning is the one of sampling from a high-dimensional probability
 distribution, e.g. the Gibbs measure $\mu_{\bA}$ of Eq.~\eqref{eq:FirstGibbs}.
 
 In the case of mean field spin glasses, physicists expect  sampling 
 to be tractable or not depending on inverse temperature $\beta$, and predict
 the location of the resulting phase transitions. For instance in the case
 of the SK model with vanishing expernal fiels (i.e. Eq.~\eqref{eq:FirstGibbs} 
 with $\bA\sim\GOE(n)$ and $h=0$) 
 Markov Chain Monte Carlo is predicted to mix rapidly for $\beta<1$
 while it is expected to fail for $\beta>1$. The critical value $\beta_c=1$,
 and similar physics predictions for other models are `deduced from'
 Parisi's formula, which also here plays an important role. 
 
 There has been recent progress towards establishing the physics 
 predictions concerning the complexity of sampling,
  but these developments are beyond the scope of this review.
 We refer to \cite{bauerschmidt2019very,eldan2022spectral,el2022sampling,celentano2024sudakov,
 huang2024sampling,kunisky2024optimality,el2025sampling} for a few pointers.

\section{Parisi's formula}
\label{sec:Parisi}

One of the landmark achievement in mathematical spin glass theory 
was the proof of Parisi's formula, which provides a variational characterization 
of the  asymptotics of $\OPT_n(\bA)$ of Eq.~\eqref{eq:OPT-Def} when
 $\bA\sim\GOE(n)$ (SK model).
In order to state Parisi's formula, we define the following class of functions:
\begin{equation}
  \cuU :=
   \Big\{ \gamma:[0,1]\to\reals_{\ge 0} \ \text{non-decreasing, right-continuous, with } \int_0^1 \!\gamma(t)\,\de t<\infty \Big\}\, .
  \label{eq:U-class}
\end{equation}
For $\gamma\in\cuU$, and $\xi:[0,1]\to \reals_{\ge 0}$
a polynomial with non-negative coefficients, let $\Phi_{\gamma,\xi}:[0,1]\times\reals\to\reals$ solve
\begin{equation}
\left\{
\begin{aligned}
  &\partial_t \Phi + \frac{1}{2}\,\xi''(t)\,\Big((\partial_x\Phi)^2 + \gamma(t)\,\partial_{xx}\Phi\Big)=0,\\
  &\Phi(1,x)=|x|.
\end{aligned}
\right.
\label{eq:parisi-pde}
\end{equation}
Define the Parisi functional
\begin{equation}
  \sP(\gamma;\xi) := \Phi_{\gamma,\xi}(0,0)\;-\;\frac{1}{2}\int_0^1 t\,\xi''(t)\,\gamma(t)\,\de t.
  \label{eq:parisi-functional}
\end{equation}

We then have the following.
\begin{theorem}\label{thm:parisi}
Let $\sP_{\cuU}^*(\xi):=\inf_{\gamma\in\cuU} \sP(\gamma;\xi)$. Then, for $H_{\sSK}(\bsigma)$ 
the Hamiltonian of Eq.~\eqref{eq:Hsigma} with $\bA\sim\GOE(n)$,
 $\xi_{\sSK}(t)=t^2/2$,   $\Sigma_n=\{+1,-1\}^n$, we have (almost surely)
\begin{equation} 
  \lim_{n\to\infty}\frac{1}{n}\,\max_{\bsigma\in\Sigma_n} H_{\sSK}(\bsigma) \;=\;
  \sP^*_{\cuU} (\xi_{\sSK})\, .
  \label{eq:opt-limit}
\end{equation}
\end{theorem}
\begin{remark}
Concentration of $\max_{\bsigma\in\Sigma_n} H_{\sSK}(\bsigma)$ around its expectation
is an elementary consequence of Borell-TIS inequality. Existence
of the limit of the expectation $\E\,\OPT_n(\bA)=n^{-1}\E\max_{\bsigma\in\Sigma_m} H_{\sSK}(\bsigma)$
is much more difficult and was proven by Guerra and Toninelli in \cite{guerra2002thermodynamic}. 
Both concentration and existence of the limit apply to the free energy density
\begin{align}
\free_n(\beta) = \frac{1}{n}\log \Big\{\sum_{\bsigma\in\{+1,-1\}^n}e^{\beta H(\bsigma)}
\Big\}\,.
\label{eq:FreeEnergy}
\end{align}
The maximum can be extracted from this quantity using the (uniform in $n,\bA$)
bound $|\OPT_n(\bA)-\free_n(\beta)/\beta|\le (\log 2)/\beta$.

Building on earlier ideas of Guerra \cite{guerra2003broken},  MT was the first to prove 
(the version of) 
Theorem \ref{thm:parisi}  for the free energy in \cite{talagrand2006parisi}. 
This proof was improved by Panchenko in \cite{panchenko2013parisi,panchenko2013sherrington}.
Auffinger and Chen \cite{auffinger2017parisi} established Theorem \label{thm:parisi} by
taking the $\beta\to\infty$ limit of the Parisi functional in \cite{auffinger2017parisi}. 
\end{remark}

\begin{remark}
  Parisi's functional turns out to be strictly convex in $\gamma$ 
  \cite{auffinger2015parisi,jagannath2016dynamic}, and hence admits a unique minimizer $\gamma_*\in \cuU$.
  Further it can be evaluated numerically with good precision
  \cite{crisanti2002analysis,schmidt2008replica,alaoui2020algorithmic}.
  Quoting the numerical estimate result of \cite{alaoui2020algorithmic}
 \begin{align}
  \sP^*_{\cuU} (\xi_{\sSK}) = 0.763168 \pm 0.000002\, .\label{eq:ParisiNumerical}
 \end{align}
 \end{remark}
 
 \begin{remark}
 Theorem \ref{thm:parisi} was stated for $\GOE(n)$ but enjoys broad 
 universality \cite{carmona2006universality,dembo2017extremal}.
 In particular consider the case of a sparse Erd\H{o}s-Renyi random graph,
 or equivalently $G\sim \SBM(n,\ovp_n,\ovp_n)$ (a stochastic block model with equal edge 
 probabilities across and within communities). Then \cite{dembo2017extremal} proves that
 \begin{align}
 \begin{cases}
  \Sigma_n:=\{\bx\in \{+1,-1\}^n:\; \<\bx,\bfone\>=0\},&\\
  \ovp_n n\to\infty&
  \end{cases} \;\;\; \Rightarrow\;\;\;
  \lim_{n\to\infty}\frac{1}{2n}\,\max_{\bsigma\in\Sigma_n} \<\bsigma,\bA_G\bsigma\> \;=\;
  \sP^*_{\cuU} (\xi_{\sSK})\, .
  \end{align}
  Following the discussion in the previous section, 
  we can use Theorem \ref{thm:parisi} to detect a latent community structure, by the following procedure.
  Given a graph $G$, we compute $\OPT_n(\bA;\Sigma_n)= \max_{\bsigma\in\Sigma_n}
   \<\bsigma,\bA_G\bsigma\>/2n$. Then we return
     \begin{align}
     T(G) = \begin{cases}
     1 & \;\;\mbox{ if }\OPT_n(\bA;\Sigma_n)\ge \sP^*_{\cuU} (\xi_{\sSK})+\delta_n,\\
     0 & \;\;\mbox{ otherwise. }
     \end{cases}
     \end{align}
     Here $T(G)=1$ means that we reject the null hypothesis of no community structure.
     Of course, this approach runs into two difficulties: $(i)$~We do not know how to set 
     $\delta_n$ (from a statistics perspective, Theorem \ref{thm:parisi} merely says that any 
     positive constant will work for $n$ large enough); $(ii)$~In general, we do not know how to 
     evaluate $\OPT_n(\bA;\Sigma_n)$. Some of the developments 
     discussed in Section \ref{sec:Algo}
     address the last problem.
 \end{remark}

\section{Computer science approaches}
\label{sec:CS-Stats}

In order to appreciate the relevance of spin glass-inspired ideas,
it is useful to highlight some of the approaches developed in 
theoretical computer science and theoretical statistics to study high-dimensional 
optimization and estimation. 
I will also highlight results that were in fact stimulated by spin glass ideas.

Again, let us focus on problem \eqref{eq:Hsigma} with $\Sigma_n=\{+1,-1\}^n$. 
It was proven in \cite{arora2005non} that this problem is hard to approximate in the worst case.
More precisely there exists a constant $c>0$ such that, unless P$=$NP, there is no 
polynomial time algorithm which is guaranteed to compute $\bsigma^{\salg} = \bsigma^{\salg}(\bA)$  such that
$H(\bsigma^{\salg})/n\ge \OPT(\bA)/(\log n)^c$.  As a consequence, no algorithm can offer the 
guarantee \eqref{eq:ApproxOPT} unless $1-\eps \le (\log n)^{-c}$. 

By and large the dominant approach towards constructing approximation algorithms
in theoretical computer science is based on convex relaxations. We begin by noting that
problem \eqref{eq:Hsigma} can be rewritten as
\begin{align}
\mbox{maximize}&\;\;\;\; \frac{1}{2} \<\bA,\bC\> \,,\label{eq:Hsigma-rel}\\
\mbox{subject to}&\;\;\; \;\bC \in \cuC_n := \conv\big\{\bsigma\bsigma^{\sT}: \;\; 
\bsigma\in\{+1,-1\}^n\big\}\, .\nonumber
\end{align} 
(In the first line  $\<\bA,\bC\> = \Tr(\bA\bC^{\sT})$ is the usual inner 
product between matrices;
in the second line $\conv(S)$ denotes the closed convex hull of the set $S\subseteq \reals^n$.)
Notice that this is equivalent to \eqref{eq:Hsigma} because the 
maximum in \eqref{eq:Hsigma-rel} is always
achieved at an extremal point $\bC=\bsigma\bsigma^{\sT}$, $\bsigma\in\{+1,-1\}^n$. The polytope $\cuC_n$
can also be characterized as the set of all second moment matrices of probability distributions over $\{+1,-1\}^n$:
\begin{align}
\cuC_n:= \Big\{\bC = \sum_{\bsigma\in \{+1,-1\}^n}\!\!\nu(\bsigma) \, \bsigma\bsigma^{\sT}:\;
\nu(\bsigma)\ge 0, \sum_{\bsigma\in \{+1,-1\}^n}\!\!\nu(\bsigma)=1\Big\}\, .\label{eq:CutPolytope}
\end{align}
Problem \eqref{eq:Hsigma-rel} is convex but intractable because checking whether $\bX\in \cuC_n$
is in itself hard. Hence, computer scientists have constructed increasingly tight relaxations,
i.e. sequences of convex sets that contain $\cuC_n$, but are computationally tractable. 
The best known example is the
elliptope $\cuE_n$ which we can define in two equivalent ways  (here
$\S^{m-1}$ denotes the unit sphere in $\reals^m$): 
\begin{align}
\cuE_n&:= \big\{\bC\in \reals^{n\times n}:\; \bC\succeq \bzero, \; C_{ii}=1\, \forall i\le n\big\}\\
& = \big\{\bC = (\<\bs_i,\bs_j\>)_{i,j\le n}\in \reals^{n\times n}:\; \bs_i \in \S^{n-1}\big\}\,, \nonumber
\end{align}
Maximizing $\<\bA,\bC\>$ over $\bC\in\cuE_n$ can be performed efficiently and
with high-degree of accuracy (using semi-definite programming \cite{boyd2004convex}). 
Once a near optimum $\bC_{\sSDP}$ is obtained, a candidate solution $\bsigma^{\sSDP}\in\{+1,-1\}^n$
can be extracted by a suitable rounding method\footnote{One option is to sample $\bg\sim \normal(0,\bC_{\sSDP})$
and let $\bsigma^{\sSDP}=\sign(\bg)$.}. A generalization of an inequality due to Grothendieck 
\cite{grothendieck1996resume,alon2006quadratic} implies that this solution satisfies
\begin{align}
\label{eq:SDPapprox}
\<\bsigma^{\sSDP},\bA\bsigma^{\sSDP}\> \ge \frac{C}{\log n}\cdot \OPT_n(\bA) \;\;\;\;\; \forall\bA\in\reals^{n\times n}\, .
\end{align}

The last result allows us to emphasize several difference between the
classical theoretical computer science viewpoint on the optimization problem \eqref{eq:Hsigma}
and the probabilistic viewpoint:
\begin{itemize}
\item[$(i)$] The most obvious difference is that the instance $\bA$ is adversarial
in computer science while it is random in spin glass theory.
A byproduct of the latter 
choice is that (under commonly chosen distributions for $\bA$), $\OPT_n(\bA)$ concentrates 
around its expectation and 
hence the problem of approximating the value $\OPT_n(\bA)$ is uninteresting.

On the other hand, finding an approximate optimizer $\bsigma^{\salg}$
as per Eq.~\eqref{eq:ApproxOPT} remains non-trivial.
\item[$(ii)$] Tools from physics and mathematical physics are tailored to 
compute structural properties such as the asymptotics of the maximum (as
in Theorem \ref{thm:parisi}), or of the free energy density \eqref{eq:FreeEnergy}.
In computer science, the focus is instead on algorithms.
\item[$(iii)$] Convex relaxations solve a problem that is --a priori-- more complex than
finding a near optimizer $\bsigma^{\salg}$.
They also provide an upper bound on the maximum value that is true deterministically. 
For instance, letting $\SDP_n(\bA) = \<\bA,\bC^{\sSDP}\>/2n$, we have that $\OPT_n(\bA)\le \SDP_n(\bA)$
for any matrix $\bA$. Together with $\<\bsigma^{\sSDP},\bA\bsigma^{\sSDP}\>\le \OPT_n(\bA)$,
we thus obtain a two-sided bound on the maximum value.

An algorithm with this property is said to provide
a `certification' or a `refutation' (since it certifies an upper bound or refutes any
lower bound contradicting it).
\end{itemize}

It is natural to ask whether the poor approximation ratio
achieved by the SDP approach (cf. Eq.~\eqref{eq:SDPapprox}) is due to the fact that 
$\bA$ is chosen adversarially. It could be that SDP works better for random instances.
As an example consider the SK  model, i.e. $\bA\sim \GOE(n)$. Then 
we have that \cite{montanari2016semidefinite} 
\begin{align}
\bA\sim \GOE(n) \;\Rightarrow\; \left\{
\begin{aligned}
\SDP_n(\bA) &= 1+o_n(1) \, ,\\ 
\frac{1}{2n}\<\bsigma^{\sSDP},\bA\bsigma^{\sSDP}\> &= 
\frac{2}{\pi}+o_n(1)\, .
\end{aligned}\right.\label{eq:SDP-random}
\end{align}
Comparing with the numerical evaluation of Parisi's formula, cf. Eq.~\eqref{eq:ParisiNumerical},
SDP achieves an approximation ratio of $2/(\pi \sP^*_{\cuU} (\xi_{\sSK}))\approx 0.834$,
which is much better than the worst case ratio of Eq.~\eqref{eq:SDPapprox}.
This leaves open the question as to whether an approximation ratio arbitrarily close to one 
can be achieved for the SK model.
 
As we briefly mentioned, convex relaxation can be systematically improved. 
Informally speaking, the elliptope $\cuE_n$ encodes all constraints on 
the cut polytope $\cuC_n$ of Eq.~\eqref{eq:CutPolytope} that can be written exclusively in terms of 
pairwise expectations of the measure $\nu$. The sum-of-squares hierarchy 
\cite{parrilo2003semidefinite,barak2014sum} yields
a sequence of convex relaxations with rich algebraic structure which eventualy
capture the complexity of $\cuC_n$, i.e.
 $\cuE_n=\cuS^{(2)}_n\supseteq  \cuS^{(4)}_n\supset \dots\supset
\cuS^{(n)}_n = \cuC_n$.

Optimization over $\cuS^{(k)}_{n}$ can be performed in time $n^{Ck}$ for some universal constant
$C$, and hence one might hope that passing from $k=2$ to some fixed $k\ge 4$
might yield an improvement over \eqref{eq:SDP-random}. Remarkably  
 \cite{kunisky2021tight,kunisky2020positivity} proved that passing from $k=2$
 to $k=n^{c}$ (for some small constant $c$) yields no improvement, i.e.
 \begin{align}
 \bA\sim \GOE(n), \; k \le n^c 
 \;\Rightarrow\; \frac{1}{2n} \max_{\bC\in\cuS^{(k)}_n}\<\bA,\bC\> = 1+o_n(1)\, .
\end{align}
Roughly speaking, even imposing all constraints on $\nu$ (cf.~Eq.~\eqref{eq:CutPolytope})
that can be written in terms of degree-$k$ moments, yields no improvement over simple degree-$2$
constraints.

Is this an indictment of convex optimization, one of most powerful tools in the hands
of the algorithm designer? We will step a bit ahead of things and anticipate that the 
answer is no. Indeed, \cite{ivkov2024semidefinite}
constructs a convex optimization-based algorithm that achieves the same approximation
ratio as the message passing  approaches discussed in the next sections
 (see also \cite{banks2021local} for an earlier related idea).
In particular, the algorithm of  \cite{ivkov2024semidefinite} yields a $(1-\eps)$ approximation
of the SK optimum, for any $\eps$, under a widely believed `no-overlap gap' conjecture
on the optimizer of Parisi's formula of Theorem \ref{thm:parisi}, discussed in Section \ref{sec:Algo}.

On the flip-side 
the algorithm of \cite{ivkov2024semidefinite} does not provide a 
`certification' (unlike relaxations of the form 
$\max_{\bC\in\cuS}\<\bA,\bC\>$ for $\cuS$ a tractable convex set).
Instead \cite{bandeira2020computational} provides evidence that any polynomial time
algorithm that takes as input $\bA$ and outputs an upper bound 
$\UB_n(\bA)\ge \OPT_n(\bA)$ that holds for every $\bA$, must be such that
$\UB_n(\bA) = 1-o_n(1)$ for $\bA\sim\GOE(n)$. In statistical 
physics language, any polynomial time algorithm 
that constructs a lower bound on the ground state energy that is always correct, must 
behave poorly for the SK model (recall that `ground state energy' 
is the minimum value of the Hamiltonian).

In order to argue that this is the case, the authors of  \cite{bandeira2020computational}
construct a distribution $\GOE_{\pm}(n)$ over matrices such that,
for  $\bA_{\pm}\sim \GOE_{\pm}(n)$ we have that $\OPT_n(\bA_{\pm}) = 1-o_n(1)$,
but a broad class of powerful algorithm fails in distinguishing between 
the two distributions $\GOE(n)$ and $\GOE_{\pm}(n)$.

\section{Message passing, belief propagation, and the cavity method}
\label{sec:MP}

The cavity method was developed in the 1980s 
\cite{mezard1986sk} as an
alternative approach to the `replica method' for deriving the the limiting properties
of the SK Gibbs measure \eqref{eq:FirstGibbs} (in particular the limiting free energy density).
The basic idea is to compare the original Gibbs measure for $n$ spins,
\begin{align}
\mu_n(\bsigma) = \frac{1}{Z_n} e^{-\beta H_n(\bsigma)}\, ,
\end{align}
to a Gibbs measure for a model in which one of the spins 
(say spin $\sigma_n$) has been removed.
In other words, we compare the Gibbs measures
$\mu_{n-1}$, $\mu_n$, for the 
two Hamiltonians 
\begin{align}
H_{n-1}(\bsigma) &:=\sum_{i<j\le n-1}A_{ij}\sigma_i\sigma_j\, ,\\
H_n(\bsigma) &:=\sum_{i<j\le n}A_{ij}\sigma_i\sigma_j=
H_{n-1}(\bsigma)+\sigma_n\sum_{i=1}^{n-1}A_{in}\sigma_i \, .
\end{align}
M\'ezard, Parisi, and Virasoro made certain assumptions about the
structure of the measure $\mu_{n-1}$,
which were inspired by the earlier replica treatment of the problem 
\cite{parisi1979infinite,parisi1980sequence}. 
These assumptions allowed them
to derive the marginal distribution of $\sigma_n$ under $\mu_n$  in terms
of marginal distributions of $\{\sigma_i\}_{i\le n-1}$ under $\mu_{n-1}$.
By postulating the existence of a large $n$ limit for certain quantities
(such as the free energy density), they were able to derive a hierarchy
of approximations for the free energy density that converged to Parisi's formula
of Theorem \ref{thm:parisi}. The $k$-th level in this hierarchy is 
commonly referred
to as `$k$-step replica symmetry breaking' 
(with the $0$-th level also known as `replica 
symmetric'.)

MT played a crucial role in translating the cavity method in a rigorous tool to 
attack spin glass models \cite{TalagrandVolI}.

Unbeknownst to physicists,  the so-called `replica symmetric cavity method'
had an algorithmic relative
known as sum-product algorithm or `belief propagation'.
This had been introduced for the first time in the 1960s by Gallager \cite{gallager1962low} as
an algorithm for iterative decoding, and subsequently rediscovered in the 1980s
in the context of graphical models' inference \cite{pearl1986fusion,koller2009probabilistic}. 

Given a finite graph $G=(V,E)$, with vertex set $V$ and edge set $E$, 
a finite alphabet $\cX$ and functions $\psi_{ij}:\cX\times \cX\to \reals_{\ge 0}$,
an `undirected graphical model' (or Markov random field)
 is the following probability 
 distribution\footnote{Belief propagation can be defined for graphical models
 with more complex structure, e.g. factor graphs or Bayesian networks, 
 but we will focus on the case of undirected graphical models for simplicity.} 
 over $\cX^V$
(well defined provided there is at least one $\bx\in\cX^V$ such that 
$\psi_{ij}(x_i,x_j)>0$ for all $(i,j)\in E$):
\begin{align}
\mu(\bx) = \frac{1}{Z}\prod_{(i,j)\in E}\psi_{ij}(x_i,x_j)\, .
\end{align}
This is of course a Gibbs measure for the Hamiltonian 
$H(\bx) = \sum_{(i,j)\in E}\log\psi_{ij}(x_i,x_j)$.

Belief propagation is an iterative algorithm to approximate low-dimensional marginals
of the graphical model $\mu$. The iteration operates on `messages' 
$\bnu = (\nu_{i\to j}^{(t)},\nu_{j\to i}^{(t)}: (i,j) \in E)$ at discrete 
times $t\in\{0,1,2,\dots \}$. Each edge $(i,j)\in E$ is given two directions 
$i\to j$ and $j\to i$ and to each direction we associate a message 
$\nu_{i\to j}^{(t)}$ which is a probability distribution over $\cX$.
Message $\nu_{i\to j}^{(t)}$
can be thought as an approximation of the marginal of $x_i$ in a model in 
which edge $(i,j)$ has been removed (call the corresponding measure joint
distribution $\mu_{(i,j)}$). 

The update equation for messages 
reads:
\begin{align}\label{eq:BP_Recursion}
\nu^{(t+1)}_{i\to j}(x_i) = \frac{1}{Z_{i\to j}}\prod_{k\in \partial i\setminus j}
\left\{\sum_{x_{k}\in \cX} \psi_{ik}(x_i,x_k)\,\nu_{k\to i}^{(t)}(x_k)\right\}\, .
\end{align}
A justification for this update rule is that, if $G$ is a tree, 
then the algorithm
computes indeed what it is supposed to compute, i.e. 
 $\nu^{(t)}_{i\to j}(\, \cdot\, ) = \mu_{(i,j)}(x_i\in \,\cdot\, )$. 
 
 Recursions analogous  to Eq.~\eqref{eq:BP_Recursion} had been known 
 in statistical physics 
  for a long time\footnote{Before the replica symmetric cavity method of 
  \cite{mezard1986sk}
   Bethe-Peierls method \cite{bethe1935statistical,peierls1936statistical} 
   had been used to study simpler 
   spin models on trees (with non-random Hamiltonians.)}. 
   However, in physics these methods  were never regarded
   algorithms (i.e. computational procedures that take as input a graphical model and
   return some quantity of interest), and instead used as analytical techniques
   (hence the  focus was on computing some average properties.)

A convenient abstraction of some important properties of belief propagation
is given by message passing algorithms, in which we retain the structure 
of the update but not necessarily the update formula \eqref{eq:BP_Recursion}.
The general message passing recursion reads:
\begin{align}\label{eq:MP_Recursion}
\nu^{(t+1)}_{i\to j}= \sF_{i\to j}\big(
\nu_{k\to i}^{(t)} :\; k\in \partial i\setminus j\big)
\, ,
\end{align}
where messages $\nu^{(t)}_{i\to j}$ are not necessarily probability distributions,
and the update functions $\sF_{i\to j}\big(\, \cdot\, )$ are arbitrary
(possibly, subject to regularity conditions.)

Between 2000 and 2005, these lines of research (cavity method in physics and 
message passing in computer science and information theory) came together, and 
each community became aware of connections between parallel developments.
In particular,
it became clear that message passing algorithms provided a rather straightforward 
path to construct algorithms motivated by physics ideas. 

An important early success of this line of work was the `survey propagation'
algorithm of \cite{mezard2002analytic}, which could be used to solve large 
random $K$-satisfiability problems and was motivated by the one-step replica-symmetry 
breaking (1RSB) cavity method. While significant amount of mathematical work was motivated 
by \cite{mezard2002analytic}, a precise understanding of the advantages and limits of survey 
propagation is still missing (see 
\cite{coja2011belief,hetterich2016analysing} 
for pointers to this literature).

\section{Approximate Message Passing and high-dimensional statistics}
\label{sec:AMP}

Message passing algorithms take as input a graph (or, equivalently, its adjacency matrix)
and --potentially-- the sets of weigth functions $\psi_{ij}$.
Approximate message passing (AMP) is the name of a class of algorithms that 
can be derived by adapting the message passing approach (which is most 
naturally suited to sparse graphs and matrices) to problems for which the input consist
of a dense matrix or dense tensor $\bA$.

One of the important properties of AMP algorithms
is that they admit a relatively simple and exact characterization
in the high-dimensional limit.

\subsection{A general algorithm and its characterization}

For the sake of simplicity, we will focus on the case in which the input is a random matrix
$\bA\sim \GOE(n)$ (see Remark \ref{rem:StateEvolution-generalizations} for extensions
to other random matrix distributions).
A specific AMP algorithm is defined by a sequence of Lipschitz functions 
$f_k:\reals^{k+1}\times\reals\to \reals$.
We define the (entrywise) action  of $f_k$
on vectors $\bx^0,\dots,\bx^k\in\reals^n$ and $\bz\in\reals^n$
by 
\begin{equation}
 f_k(\bx^0,\dots,\bx^k;\bz) :=
  \begin{pmatrix}
    f_k(x^0_1,\ldots,x^k_1; z_1)\\
    \vdots\\
    f_k(x^0_n,\ldots,x^k_n; z_n)
  \end{pmatrix}.
  \label{eq:separable-fk}
\end{equation}
We also write $(\bx^0,\ldots,\bx^k)=\bx^{\le k}$ and $f_k(\bx^0,\ldots,\bx^k;\bz)=f_k(\bx^{\le k};\bz)$.

\begin{definition}[AMP iteration; Matrix data case]\label{def:AMP}
The \emph{AMP iteration} with nonlinearities $f_k$ is the sequence $(\bx^k)_{k\ge 0}\subset\reals^n$ defined by
\begin{equation}
  \bx^{k+1} = \bA\, f_k(\bx^{\le k}; \bz)\;-\; \sum_{j=0}^{k} \hsb_{k,j}\, f_j(\bx^{\le j}; z)\, ,
  \label{eq:amp-iter}
\end{equation}
where the \emph{Onsager coefficients} are estimated by
\begin{equation}
 \hsb_{k,j} :=\frac{1}{n}\sum_{i=1}^n \partial_j f_k(x^{\le k}_i; z_i),
  \qquad j=0,\ldots,k,
  \label{eq:bhats}
\end{equation}
and $\partial_j f_k$ denotes the weak derivative of $f_k$ with respect to 
$x^{j}_i$.
\end{definition}

In what follows we will assume that $\bx^0,\bz$ are either deterministic or random and
independent of $\bA$, and that the empirical distribution of $((x^0_i,z_i):\, i\le n)$ converges
in $W_2$ (Wasserstein-$2$) distance to a law $p_{X_0,Z}$ with finite second moment
\begin{align}\label{eq:InitW2}
\frac{1}{n}\sum_{i=1}^n \delta_{x^0_i,z_i} \stackrel{W_2}{\Rightarrow} p_{X_0,Z}\, .
\end{align}
(If $\bx^0$, $\bz$ are random, the above convergence is assumed to hold almost surely.)
This convergence is also more concretely stated in terms of test functions.
We say that $\psi:\reals^{m}\to\reals$  is pseudo-Lipschitz of order two 
(and write $\psi\in\PL$) if 
$|\psi(\bx)-\psi(\by)|\le C\big(1+\|\bx\|_2+\|\by\|_2\big)\|\bx-\by\|_2$
for all $\bx, \by\in\reals^m$. 
Then Eq.~\eqref{eq:InitW2} is equivalent \cite{villani2008optimal} to the statement that, 
for any $\psi\in\PL$,
\begin{align}
\lim_{n\to\infty} \frac{1}{n}\sum_{i=1}^n \psi(x^0_i,z_i) = \E\{\psi(X_0,Z)\}\, .
\end{align}

\begin{definition}[State Evolution process and recursion]\label{def:SE}
The \emph{state evolution} process is the stochastic process $((X_k)_{k\ge 0},Z)$
where $(X_k)_{k\ge 1}$ is Gaussian and independent of $(Z,X_0)\sim p_Z\times p_{X_0}$, with 
covariance defined recursively by (here $X^{\le k}=(X_0,X_1,\dots,X_k)$)
\begin{equation}
  Q_{k+1,j+1} \;=\; \E\Big\{ f_k(X^{\le k}; Z)\, f_j(X^{\le j}; Z)\Big\}\, .
  \label{eq:SE-Q}
\end{equation}
\end{definition}

The following theorem gives an asymptotic characterization 
of the AMP iteration. It was  proven in this form in \cite{BM-MPCS-2011,javanmard2013state},
generalizing a Gaussian conditioning approach from \cite{bolthausen2014iterative}.
\begin{theorem}\label{thm:StateEvolution}
  Assume Eq.~\eqref{eq:InitW2} holds and that the functions $f_k$  are Lipschitz continuous.  
Then, for any $\psi\in\PL$, almost surely,
\begin{equation}
  \lim_{n\to\infty}\frac{1}{n}\sum_{i=1}^n \psi\big(x^{\le k}_i, z_i\big)
  \;=\; \E\,\psi\big(X^{\le k}, Z\big)\, .
  \label{eq:lln}
\end{equation}
\end{theorem}
Informally, Theorem \ref{thm:StateEvolution} establishes that the 
random variables $Z,X_0,X_1,\dots,X_k$ are the infinite-dimensional 
equivalent of vectors $\bz,\bx^0,\bx^1,\dots,\bx^k$. 
In particular, the vectors
$\bz,\bx^0,\bx^1,\dots,\bx^k\reals^n$ (endowed with the 
normalized scalar product $\<\cdot,\cdot\>/n$)
are approximately isometric to 
$Z,X_0,X_1,\dots,X_k$ (endowed with the $L^2$ scalar product).
For instance, for any fixed $i,j$,
\begin{align}
\lim_{n\to\infty}\frac{1}{n}\<\bx^i,\bx^j\> = \<X_i,X_j\>_{L^2}\, .
\end{align}

\begin{remark}\label{rem:StateEvolution-generalizations}
Several generalized/stronger forms of this theorem have been proved in recent years.
 Versions of AMP exist for rectangular matrices with independent entries of
unequal variances  \cite{javanmard2013state},
tensors \cite{alaoui2020optimization}, non-separable functions $f_k$ (i.e. functions that do not act coordinatewise)
\cite{berthier2019state}, multi-matrix iterative algorithms \cite{gerbelot2023graph}.

State evolution has been proved for matrices with non-Gaussian independent entries
 \cite{bayati2015universality,chen2021universality};
distributed according to orthogonally invariant ensembles
 \cite{fan2022approximate,wang2024universality}
(at the condition of modifying the definition of coefficients $\hsb_{k,j}$);
and even semirandom matrices 
\cite{dudeja2023universality,wang2024universality}.
\end{remark}

\begin{remark}
The SE analysis can be strengthened to obtain quantitative tail bounds and 
control over $O(n^{\delta})$ iterations for some $\delta>0$ \cite{rush2018finite,li2023approximate,li2022non}. 
\end{remark}

\subsection{Low-rank matrix estimation via AMP}

\begin{figure}[t]
  \includegraphics[width=0.5\linewidth]{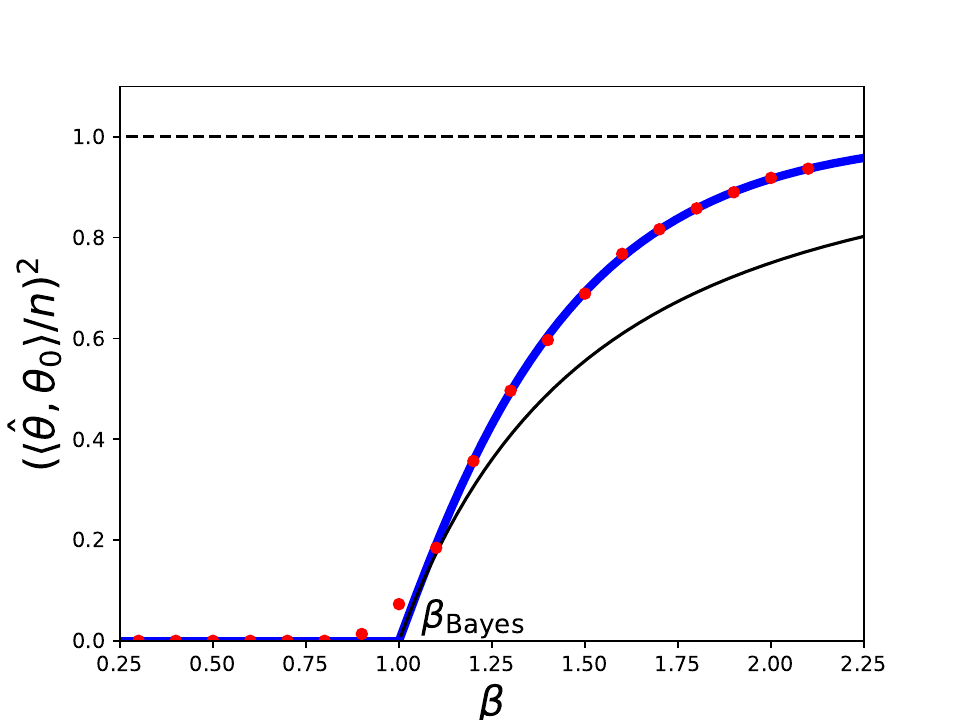}
  \includegraphics[width=0.5\linewidth]{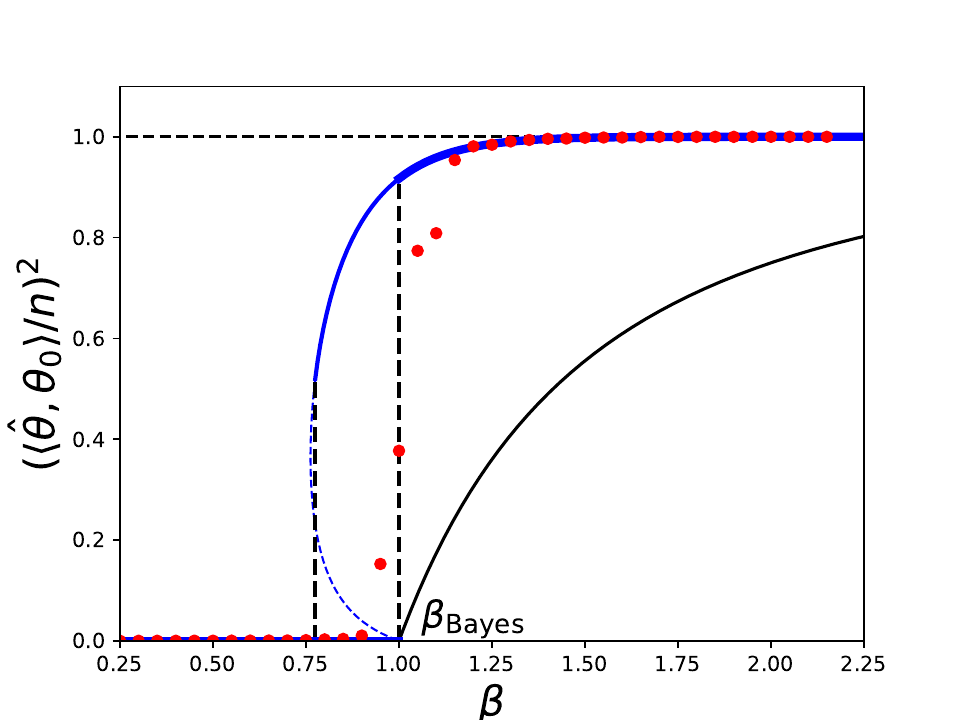}
  
  \caption{}\label{eq:LowRank}
 \end{figure}

To demonstrate the application of AMP and spin glass
ideas to high-dimensional statistics, it is instructive to consider a
simple rank-one estimation problem. We want to estimate a vector $\btheta_*\in\reals^n$
from a single noisy observation $\bY\in\reals^{n\times n}$ of the outer product $\btheta_*\btheta_*^{\sT}$:
\begin{equation}
  \bY \;=\; \frac{\lambda}{n}\,\btheta_* \btheta_*^\sT + \bA,\qquad \bA\sim\GOE(n).
  \label{eq:spiked-model}
\end{equation}
We observe $\bY$ and our goal is to estimate $\btheta_\ast$. We will assume the coordinates 
of $\btheta_*$ to have empirical distribution  converging in $W_2$
\begin{align}
\frac{1}{n}\sum_{i=1}^n\delta_{\theta_{*, i}} \stackrel{W_2}{\Rightarrow}
 P_{\Theta_\ast}\, .
 \end{align}
 We will assume $P_{\Theta_*}$ to have unit second moment $\E\{\Theta^2_\ast\}=1$,
 (here and below $\Theta_*\sim  P_{\Theta_*}$). A special example is the one in which 
the coordinates of $\btheta_*$
 are themselves random $(\theta_{*,i}:i\le n)\sim_{iid}P_{\Theta_*}$.

One of the motivations for considering this specific
problem resides in the fact that the connection to standard 
spin glass models is particularly straightforward. Indeed, under
the i.i.d. prior $\btheta\sim P_{\Theta}^{\otimes n}$, 
the posterior distribution of $\btheta$ given $\bY$ is:
\begin{align}
  P_{\btheta|\bY}(\de\btheta|\bY) = \frac{1}{Z_n(\bY)}
  \exp\Big\{\frac{\lambda}{2}\<\btheta,\bY\btheta\>-\frac{\lambda^2}{4n}\|\btheta\|_2^4
  \Big\}
P_{\Theta}^{\otimes n}(\de\btheta)\, .\label{eq:FirstPosterior}
\end{align}
This is very similar to the Gibbs measure of a spin model,
e.g. the SK model given by \eqref{eq:FirstGibbs}
whereby $H(\bsigma)=\<\bsigma,\bA\bsigma\>/2$ for $\bA\sim\GOE(n)$.
The main differences are:
$(i)$~The different reference measure $P_{\Theta}^{\otimes n}(\de\btheta)$;
$(ii)$~The different distribution of the random matrix $\bY$, which
is given by Eq.~\eqref{eq:spiked-model}.

A general AMP algorithm produces a sequence of estimates $(\btheta^k)_{k\ge 1}$
via
\begin{align}
  \btheta^{k+1} &= \bY f_k\big(\btheta^{\le k}\big) - 
  \sum_{j=1}^{k} \hsb_{kj}f_{j-1}(\btheta^{\le j-1})\,. \label{eq:g-raw}
\end{align}
with  $\hsb_{kj}$ defined as in Eq.~\eqref{eq:bhats} and initialization $\btheta^0 = \E\{\Theta_*\}\bfone$
(this is the best estimate of  $\btheta_*$ absent any  information other than the prior $P_{\Theta_*}$.)
Substituting the form of $\bY$ of Eq.~\eqref{eq:spiked-model}, we get
 \begin{align} 
  \btheta^{k+1} &= \frac{\lambda}{n}\,\<\btheta_*, f_k(\btheta^{\le k})\>\btheta_* \;+\; \bA\, 
  f_k(\btheta^{\le k}) \;-\; \sum_{j=0}^{k-1} \hsb_{k,j}\, f_j(\btheta^{\le j}). \label{eq:g-expanded}
\end{align}
This differs from the AMP iteration of Eq.~\eqref{eq:amp-iter} because of the shift 
$(\lambda/n)\,\<\btheta_*, f_k(\btheta^{\le k})\>\btheta_*$. However it is possible to 
show that the coefficient $\lambda\<\btheta_*, f_k(\btheta^{\le k})\>/n$ 
concentrates around its expectation (call it $\omu_k$) for any constant $k$,
and hence we can apply Theorem \ref{thm:StateEvolution} to $\bx^k = \btheta^k - \omu_k\btheta_*$.

As a consequence of state evolution, 
we can construct a sequence of random variables $(\Theta_k:k\ge 0)$ in the same probability
can construct a sequence of random variables $(\Theta_k:k\ge 0)$ in the same probability
space as $\Theta_*$ such that, for any $\psi\in\PL$,
\begin{equation}
  \lim_{n\to\infty}\frac{1}{n}\sum_{i=1}^n \psi(\theta_{*,i}, \theta_i^{\le k}) =
  \E\big\{\psi(\Theta_*, \Theta^{\le k})\big\}\, .
  \label{eq:empirical-theta}
\end{equation}
Because of the shift in Eq.~\eqref{eq:g-expanded}, the random variables $\Theta_j$
have a positive component along $\Theta_*$:
\begin{align}
  \Theta^{\le k} &= \mu^{\le k} \Theta_\ast + G^{\le k}, \;\;\;\;\;\;\Theta_*\indep G^{\le k},\;
  \; G^{\le k}\sim\normal(0,Q_{\le k,\le k})\,,\label{eq:theta-se}
 \end{align}
where $\Theta^{\le k}= (\Theta^1,\dots,\Theta^k)$, $\Theta^0=\E\{\Theta_*\}$ is non-random, 
and $\mu, Q$
are determined by the following recursion for $k,j\ge 0$
\begin{align}
  \mu_{k+1} &= \lambda\, \E\,\Big\{\Theta_*\, f_k(\Theta^{\le k})\Big\}, \label{eq:mu-se}\\
  Q_{k+1,j+1} &= \E\,\Big\{ f_k(\Theta^{\le k})\, f_j(\Theta^{\le j})\Big\}. \label{eq:Q-se}
\end{align}

A possible use of this characterization is to design the optimal nonlinearities $f_k$. 
We want maximize the correlation:
\begin{equation}\label{thm:opt-fk}
  \rho_{k,n} := \frac{\< \btheta_*, \btheta^k\>}{\|\btheta_*\|\,\|\btheta^k\|} \,.
\end{equation}
State evolution allows to compute the limit
\begin{align}
  \plim_{n\to\infty}\rho_{k,n}=\rho_k:=\frac{\mu_k}{\sqrt{\mu_k^2+Q_{k,k}}} .
\end{align}

In fact the optimal sequence of functions $f_k$ can be determined explicitly,
as stated in the following theorem from 
\cite{celentano2020estimation,montanari2022statistically}.
\begin{proposition}\label{thm:opt-fk}
For any $K\ge 0$, the asymptotic correlation $\rho_{K}$ is uniquely
maximized if we set recursively, for  $k\in \{0,\dots , K-1\}$, 
\begin{equation}
  f_k(\Theta^{\le k}) \;=\; h_k(\Theta^k), \qquad 
  h_k(y) \;=\; \E\{\Theta_*\mid \mu_k \Theta_* + \tau_k G = y\},
  \label{eq:opt-denoiser}
\end{equation}
where $G\sim\normal(0,1)$, and $\mu_k$, $\tau_k = \sqrt{Q_{kk}}$ are given recursively by
\begin{align}
  \mu_{k+1} &= \lambda\, \E\,\Big\{\Theta_\ast\, h_k\!\big(\mu_k \Theta^\ast + \tau_k G\big)\Big\}, \;\;\;\;\;
  \tau_{k+1}^2 = \E\,\Big\{ h_k\!\big(\mu_k \Theta^\ast + \tau_k G\big)^2 \Big\}\, , \label{eq:tau-opt}
\end{align}
with initialization $\mu_1=\lambda\E\{\Theta_*\}^2$, $\tau^2_1=\E\{\Theta_*\}^2$.
Equivalently, for $\gamma_k = (\mu_k/\tau_k)^2$, $\gamma_1=\lambda\E\{\Theta_*\}$, we have
\begin{equation}
  \gamma_{k+1} \;=\; \lambda^2\, \sF(\gamma_k), \qquad \sF(\gamma) := 
  \E\left\{\E\{\Theta_*|\sqrt{\gamma_k}\Theta_*+G\}^2\right\}.
  \label{eq:gamma-recursion}
\end{equation}
The corresponding correlation is $\rho_k = \gamma_k/\sqrt{1+\gamma_k^2}$, and we have
$\mu_k=\gamma_k/\lambda$, $\tau^2_k=\gamma_k/\lambda^2$,
and for any $k,j\ge 1$, $Q_{k,j} = \tau^2_{k\wedge j}$.
\end{proposition}
In words, the optimal estimation achieved by a general first order method
after any constant number of iterations is predicted by the simple one-dimensional recursion 
\eqref{eq:gamma-recursion}.

\begin{remark}[Optimality within broader classes of algorithms]
The optimal AMP algorithm is also optimal within the 
`generalized first order methods' \cite{celentano2020estimation},
as well as among `low degree polynomials' 
\cite{montanari2025equivalence}.

Informally, a generalized first order method computes recursively a sequence
$\btheta^k\in\reals^n$ where, for each $k$, 
$\btheta^{k+1} = \bY f_k(\btheta^{\le k})+g_k(\btheta^{\le k})$
for two sequence of functions $(f_k:k\ge 1)$, $(g_k:k\ge 1)$ that operate
row-wise on the matrices $[\btheta^1,\dots,\btheta^k]$.
The final estimator is $\hbtheta= \btheta^K$ for fixed 
number of iterations $K$.  A number of 
classical iterative methods fit in this framework, e.g. gradient descent 
with respect to standard regularized loss functions. 

Low degree polynomials are instead defined by the condition that
each coordinate $\htheta_i(\bY)$ is a polynomial  of maximum degree
$D$ in the entries of $\bY$. 

The maximum accuracy achivable by either generalized first order methods
(with $K=O(1)$) or low-degree polynomials (with $D=O(1)$) is the same as the 
one achieved by AMP.  
\end{remark} 

\subsection{Optimal low-rank matrix estimation and 
statistical computation gaps}

It is possible to prove that the function $\gamma\mapsto \sF(\gamma)$
is monotone non-decreasing and therefore
the iteration \eqref{eq:gamma-recursion} always converges to 
fixed points.
We define the smallest stable fixed point 
(note that, under the assumptions, $\gamma\mapsto 
\lambda^2\sF(\gamma)$ is continuous):
\begin{align}
\gamma_{\salg}(\lambda) := \inf\big\{\gamma>0:\; 
\lambda^2\sF(\gamma)<\gamma\big\}\, ,
\end{align}
and $\rho_{\salg}:= \gamma_{\salg}/\sqrt{1+\gamma_{\salg}^2}$.
In particular, if $\lambda^2\sF'(0)=\lambda^2>1$,
we have $\lim_{\gamma\to\infty}\gamma_k=\gamma_{\salg}(\lambda)$
for all small enugh initialization $\gamma_0>0$. 
In other words, $\gamma_{\salg}(\lambda)$ appears to be a barrier
that cannot be improven upon, even when AMP
is initialized with a correlation with $\btheta_*$.
As discussed below, $\rho_{\salg}$ is conjectured
to be the optimal accuracy achievable by any polynomial time estimator.

Remarkably the optimal accuracy achievable by
\emph{any} estimator (i.e. any measurable function $\bY\mapsto \hbtheta(\bY)$)
is also given by a fixed point of the map $\gamma\mapsto \lambda^2\sF(\gamma)$.
In order to state this result precisely, we define the following 
function:
\begin{equation}
  \Psi(\gamma;b,\lambda,P_{\Theta}) := \frac{1}{4} \big(q-\lambda^2\big)^2-
  \frac{1}{2} b \gamma+\Info(\gamma;P_{\Theta}) \, ,\label{eq:PsiDef}
\end{equation}
\begin{equation}
 \Info(\gamma;P_{\Theta}) := \E\log \frac{\de p_{\Yeff \vert \Theta}}{\de p_{\Yeff}}\,,
 \;\;\;
   \Yeff = \lambda\sqrt{\gamma}\Theta +G \, ,\;\;\; (\Theta,G)\sim P_{\Theta}\otimes 
   \normal(0,1)\, .
\end{equation}
We then define
\begin{equation}
  \gamma_{\sBayes}(\lambda) := \arg\min\{ \Psi(\gamma;b=0,\lambda,P_{\Theta}):\;\;\gamma>0\big\}\, ,
\end{equation}
and $\rho_{\sBayes}= \gamma_{\sBayes}/\sqrt{1+\gamma_{\sBayes}^2}$.

It is immediate to show  that $\Psi'(\gamma;0,\lambda,P_{\Theta})=0$
if and only if $\gamma=\lambda^2\sF(\gamma)$:
critical point of $\Psi(\gamma;0,\lambda,P_{\Theta})$
correspond to fixed points of $\gamma\mapsto \lambda^2\sF(\gamma)$.

We can then state some fundamental facts about optimal estimation in
the model \eqref{eq:spiked-model}
(summarizing results from 
\cite{lelarge2016fundamental,montanari2017estimation}, see 
\cite{montanari2025equivalence}).
\begin{theorem}\label{thm:SummarySpiked}
If either $\E\{\bTheta\}\neq 0$ or  $\E\{\bTheta\}= 0$
and $\lambda>1$, then 
there exists an algorithm  with complexity $O(n^2\log n)$
that computes $\hbtheta^{\salg}=\hbtheta^{\salg}(\bY)$
such that
\begin{align}
  \label{eq:AlgAchievable}
\lim_{n\to\infty} 
\frac{|\<\hbtheta^{\salg},\btheta_*\>|}{\|\btheta_*\|\|\hbtheta^{\salg}\|}
= \rho_{\salg}\, . 
\end{align}
On the other hand, if $b\mapsto 
\Psi_*(b,\lambda,P_{\Theta}):=\min_{\gamma>0}\Psi(\gamma;b,\lambda,P_{\Theta})$
is differentiable at $b=0$, then the optimal estimator
(i.e. the one maximizing
$\E\{|\<\hbtheta,\btheta_*\>|/\|\btheta_*\|\|\hbtheta\|\}$)
is  such that
\begin{align}\label{eq:StatAchievable}
  \lim_{n\to\infty} 
  \frac{|\<\hbtheta^{\sBayes},\btheta_*\>|}{\|\hbtheta_*\|\|\hbtheta^{\sBayes}\|}
  = \rho_{\sBayes}\, . 
\end{align}
\end{theorem}
\begin{remark}
  When $\E\{\Theta\}\neq 0$, we hane $\gamma_0>0$ and therefore
  $\lim_{k\to\infty}\gamma_k=\gamma_{\salg}$. In this case,
  Eq.~\eqref{eq:AlgAchievable} in 
  Theorem~\ref{thm:SummarySpiked} follows from Propositon~\ref{thm:opt-fk}.

When $\E\{\Theta\}=0$, state evolution predicts $\gamma_k=0$ 
for all fixed number of iterations $k$.
In this case, new ideas are  necessary to prove 
 claim \eqref{eq:AlgAchievable}.
This is proved using a spectral initialization 
for AMP in \cite{montanari2021estimation}, and by 
analyzing AMP beyond $O(1)$ iterations in 
\cite{li2023approximate}.
\end{remark}

\begin{remark}
The function $\Psi_*(0;\lambda,P_{\Theta})=
\min_{\gamma>0}\Psi(\gamma;0,\lambda,P_{\Theta})$ of Eq.~\eqref{eq:PsiDef}
possesses an intrinsic meaning, since it coincides (up to a simple additive
constant) with the asymptotic value of the free energy density
$n^{-1}\log Z_n(\bY)$ for the Gibbs measure of Eq.~\eqref{eq:FirstPosterior}.

Indeed, the crucial step in establishing \eqref{eq:StatAchievable}
is to characterize the limit of the free energy density $n^{-1}\log Z_n(\bY)$.
One possible approach towards this goal 
---pursued in \cite{lelarge2016fundamental}---
is to use the cavity method as developed in \cite{TalagrandVolI}.
This is particularly effective because concentration of the overlap
can be established \emph{a priori} using the Bayes 
structure\footnote{Identities that follow from this structure are 
also called `Nishimori identities' in physics 
\cite{iba1999nishimori}.} of the 
measure \eqref{eq:FirstPosterior}.
\end{remark}

%
%
\section{Algorithimc thresholds for spin glasses}
\label{sec:Algo}

Armed with the general tools of Section \ref{sec:AMP}
(AMP and its state evoluton analysis), we revisit the 
problem of optimizing spin-glass Hamiltonians. 
Namely, we want to develop an algorithm that takes as input
(a specification of) the Hamiltonian $H(\,\cdot\,)$ and returns
$\bsigma^{\salg}$ such that (cf. Eq.~\eqref{eq:ApproxOPT})
\begin{align}
\lim_{n\to\infty} \prob\Big\{H(\bsigma^{\salg})\ge (1-\eps)
\max_{\bsigma\in\Sigma_n}
H(\bsigma)\Big\} =1\, .   \label{eq:HamApproxProblem}
\end{align}
As discussed in Section \ref{sec:CS-Stats}, 
we would like $\eps$ as small as possible, with 
SDP relaxations achieving $(1-\eps)\approx 0.834$ for the SK Hamiltonian.

\subsection{Optimizing spin glass Hamiltonians}

It is useful to generalize the discussion of the optimization problem 
to the case of a mixed $p$-spin model.
This is a  natural generalization of Eq.~\eqref{eq:Hsigma}
whereby the Hamiltonian is a polynomial with Gaussian coefficients.
Namely
\begin{align}\label{eq:MixedPSpinModel}
H(\bsigma) = \sum_{k=2}^{k_{\max}}\frac{\sqrt{\xi_k}}{n^{(k-1)/2}}\<\bG^{(k)},\bsigma^{\otimes k}\>\, .
\end{align}
where $\bG^{(k)} =(G^{(k)}_{i_1\cdots i_k})_{i_1\cdots i_k\le n}\sim_{iid}\normal(0,1)$.
Equivalently, $H(\;\cdot\; )$ is a centered Gaussian process with covariance
\begin{align}
\E\big\{H(\bsigma_1)H(\bsigma_2)\big\} = n\xi\Big(\frac{\<\bsigma_1,\bsigma_2\>}{n}\Big)\, \, ,
\end{align}
where $\xi(t)=\sum_{k=2}^{k_{\max}}\xi_k t^{k}$ is a 
polynomial\footnote{The results can be easily
 extended to $k_{\max}=\infty$
under a suitable decay condition on the coefficients $\xi_k$.}.

  Theorem \ref{thm:parisi} remains verbatimly true for the Hamiltonian
\eqref{eq:MixedPSpinModel}, provided $\xi$ in Parisi's 
formula matches the polynomial $\xi$ in Eq.~\eqref{eq:MixedPSpinModel}.

Parisi's formula and its proof in 
\cite{talagrand2006parisi,panchenko2013parisi,auffinger2015parisi}
are non-constructive and do not provide a solution to problem
\eqref{eq:HamApproxProblem}. In order to state results on the latter,
we define the class of functions
\begin{equation}
  \cuL
   := \Big\{ \gamma:[0,1])\to\reals_{\ge 0} : 
   \|\xi''\gamma\|_{\sTV[0,t]}<\infty \, \forall t\in (0,1)\, ,\;\;\;
   \int_0^1\xi''\gamma(t)\,\de t<\infty \Big\}\, .
  \label{eq:L-class}
\end{equation}
We note that $\cuL$, is a superset of the set of functions $\cuU$ introduced 
in Eq.~\eqref{eq:U-class}.
The following theorem was proven in \cite{subag2021following}
for $\Sigma_n=\sqrt{n}\S^{n-1}$ and in
\cite{montanari2019optimization,alaoui2020optimization}
for  $\Sigma_n=\{+1,-1\}^n$. 
\begin{theorem}\label{thm:alg-guarantee}
Let $\Sigma_n\in\big\{\sqrt{n}\S^{n-1},\{+1,-1\}^n\big\}$,
and $\sP_{\cuL}^*(\xi):=\inf_{\gamma\in\cuL} \sP(\gamma;\xi)$,
where $\sP(\gamma;\xi)$ si defined as in Eq.~\eqref{eq:parisi-functional}
for $\Sigma_n=\{+1,-1\}^n$, and replacing the boundary condition
in Eq.~\eqref{eq:parisi-pde} by $\Phi(1,x)=x^2/2$ for
$\Sigma_n=\sqrt{n}\S^{n-1}$.

Assume the infimum $\inf_{\gamma\in\cuL} \sP(\gamma;\xi)$
is achieved for some 
$\gamma^\ast\in\cuL$. 
Then for any $\eps>0$ there exists an 
algorithm  with complexity at most $C(\eps)\, n^{k_max}$  such that
\begin{equation}
  \lim_{n\to\infty}\prob\left(
  \frac{1}{n}H(\sigma^{\mathrm{alg}}) \ge \sP_{\cuL}^*(\xi) - \eps\right) = 1\ .
  \label{eq:alg-guarantee}
\end{equation}
\end{theorem}

\begin{remark}[Spherical model]
For the case $\Sigma_{n}=\sqrt{n}\S^{n-1}$, the Parisi PDE can be solved 
explicitly ($x\mapsto  \Phi(t,x)$ is a quadratic function for all $t$.)
The minimization over $\gamma$ can be worked out explicitly
and the assumption that $\inf_{\gamma\in\cuL} \sP(\gamma;\xi)$ 
is achieved holds, resulting in the formula:
\begin{align}
\Sigma_n=\sqrt{n}\S^{n-1}\;\;\;\Rightarrow\; \;\;
 \sP_{\cuL}^*(\xi) = \int_0^1\!\sqrt{\xi''(t)}\; \de t\, .
\end{align}
\end{remark}

\begin{remark}
Since $\cuU\subseteq \cuL$, we have $\sP_{\cuU}^*(\xi)\le \sP_{\cuL}^*(\xi)$:
the algorithmically achievable value is no larger than the maximum,
a useful sanity check. 

Further, it can be shown that,
when the infimum over $\cuU$ is achieved at a $\gamma^*$ that is 
strictly increasing on $[0,1)$,  then $\gamma^*$ also achieves the 
infimum over the larger space $\cuL$, and hence 
$\sP_{\cuU}^*(\xi) = \sP_{\cuL}^*(\xi)$. 

The condition that $\gamma^*$ is strictly increasing
is also referred as `no overlap gap,' since $\gamma^*$ is related 
to the overlap distribution (a central object of study in spin glass theory) 
at low temperature. Within statistical physics this
condition is believed to hold for the SK model $\xi(t) = \xi_2t^2$.
\end{remark}

The proof of Theorem \ref{thm:alg-guarantee} in 
\cite{montanari2019optimization,alaoui2020optimization} is based on the construction
of a class of AMP algorithms with certain specific properties,
also known as incremental AMP, or IAMP. 
(We also refer to \cite{jekel2025potential} for an alternative 
approach based on the TAP free energy.)
The value these algorithms achieve can be 
characterized in terms of a stochastic optimal control problem.
Theorem \ref{thm:alg-guarantee}  then follows from strong duality 
for this control problem.
Even when strong duality does not hold, this argument
yields an achievable value.
We will state this stochastic optimal control problem
for the case $\Sigma_n=\{+1,-1\}^n$.
\begin{theorem}\label{thm:OPT-Control}
For  $(B_t:t\ge 0)$
  a standard Brownian, and let $\cuD_0$ be the 
  set of all processes that are progressively measurable
  with respect to the natural filtration of $(B_t:t\ge 0)$. 
  Further define two spaces
  of control processes, $\cuD$, $\cuD_s$
   (with $\cuD_s\subseteq \cuD\subseteq \cuD_0$) as follows:
\begin{align}
  \cuD&:=\Big\{U\in\cuD_0: \; \;
  \xi''(t)\E \{U^2_t\} = 1\;\; \forall t\in [0,1],\; 
  \int_0^1 \sqrt{\xi''(t)} \, U_t\,\de B_t\in (-1,1)\Big\}\, ,\\
  \cuD_s&:=\Big\{U\in\cuD: \;
  U_t = u(t,X_t)\, ,
   \de X_t= v(t,X_t)\de t+\de B_t,\, \forall t\in[0,1], X_0=0, 
   u,v\;\mbox{\rm Lipschitz}\Big\}\,.
\end{align}
And $\sV:\cuD\to \reals$:
\begin{align}
\sV(U) := &\int_0^1\,\xi''(t)\,\E\, U_t\,  \de t\, .
\end{align}
Then, for any $U\in\cuD_s$, $\eps>0$ there exists an algorithm with 
complexity at
 most $C(\eps)n^{k_{\max}}$, and output $\bsigma^{\salg}\in \Sigma_n$
  such that $n^{-1}H(\bsigma^{\salg})\ge \sV(U) - \eps$ with 
  probability converging to one as $n\to\infty$. 
\end{theorem}

As mentioned above, Theorem \ref{thm:alg-guarantee} can be proven from
Theorem \ref{thm:OPT-Control} via a duality argument.
Namely, Theorem \ref{thm:OPT-Control}  implies that 
$\ALG_s:=\sup_{U\in\cuD_s} \sV(U)$ is an achievable value,
while of course $\ALG_s\le \ALG:=\sup_{U\in\cuD} \sV(U)$.
Minimization of the Parisi functional $\sP(\gamma;\xi)$ turns out to be the dual of 
the  maximization problem $\sup_{U\in\cuD}\sV(U)$,
(in particular $\ALG_s\le \ALG \le \sP_{\cuL}^*(\xi)$). 
When the infimum  is achieved, it can be shown to yield 
$U_*\in\cuD_s$ such that $\sV(U_*)=\sP_{\cuL}^*(\xi)$, and
therefore $\ALG_s = \ALG =  \sP_{\cuL}^*(\xi)$.

The proof of Theorem \ref{thm:OPT-Control} proceeds
by constructing an AMP algorithm with the following iteration:
\begin{align}
  \bz^{\ell+1} &= \nabla H(\bm^{\ell}) \;-\; 
  \sum_{j=1}^{\ell} b_{\ell j}\, \bm^{j-1},
  \;\;\;\;\;\;
  \bm^{\ell} = f_\ell\!\big(\bz^{\le \ell}\big)\, .
  \label{eq:iamp-m}
\end{align}
Note that the gradient takes the form
\begin{align}
  \nabla H(\bm) = \sum_{k=2}^{k_{\max}} 
  k\sqrt{\xi_k} \bG^{(k,s)}\{\bm\}\, .
\end{align}
where $\bG^{(k,s)}$ is the symmetrization of $\bG^{(k)}$, and 
$\bG^{(k,s)}\{\bm\}$ denotes the contraction of $\bm$
with $k-1$ indices in $\bG^{(k,s)}$. Note that this is  a natural 
generalization of the basic AMP iteration of Eq.~\eqref{eq:amp-iter} 
to   tensor data.

In the 1980s physicists developed a fascinating picture of the 
geometry of near optima of spin glasses 
\cite{mezard1984nature,mezard1984replica,mezard1985microstructure}.
A number of elements of this picture were proved since MT's 
breakthrough proof of Parisi's formula \cite{talagrand2006parisi},
see for instance 
\cite{talagrand2010construction,panchenko2013parisi,chen2023generalized,chen2021generalized}.
While these results are not used in a formal way in the proof of
Theorem \ref{thm:alg-guarantee} and Theorem \ref{thm:OPT-Control},
they provide the key intuition to construct a good family of functions
 $f_k$ in Eq.~\eqref{eq:iamp-m}. It is therefore useful to 
 review this picture in some detail: we will do it at a heuristic level.

Consider the superlevel set 
\begin{align}
\cS_n(\eps):=\Big\{\bsigma\in \Sigma_n: \, H(\bsigma)\ge
(1-\eps)\max_{\bsigma\in\Sigma_n} H(\bsigma)\Big\}\, .
\end{align}
A crucial property of this set is that, with high probability,
 it contains an arbitrarily
large (as $n\to\infty$) finite set $\cS^*_n(\eps)
\subseteq \cS_n(\eps)$ which is:
$(i)$~Ultrametric, i.e. 
$\|\bsigma_1-\bsigma_2\|\le \max(\|\bsigma_1-\bsigma_3\|,\|\bsigma_2-\bsigma_3\|)$
for all $\bsigma_1,\bsigma_2,\bsigma_3\in \cS^*_n(\eps)$;
$(ii)$~Roughy balanced, i.e. the number of points of $\cS^*_n(\eps)$ in any 
ball $\Ball(\bsigma;r)$ is roughly the same for all $\bsigma\in \cS^*_n(\eps)$;
$(iii)$~Centered, i.e. the barycenter of $\cS^*_n(\eps)$ is close to the origin.

As a consequence of these properties, the near optima in
$\cS^*_n(\eps)$ can be organized as the leaves of a tree.
We construct this tree by recursively clustering the points of $\cS^*_n(\eps)$
as follows.
We partition the interval $[0,1]$ $1/\delta$ subintervals of length 
$\delta$.
At each level $q\in \{0,\delta,\dots, 1-\delta,\delta\}$ the 
construction
generates clusters $\cC_\ell$, with centers $\bm_{\ell}$, indexed by
$\ell\in L(q)$,
starting with clusters formed of single points at level $q=1$.
To construct level-$q$ clusters, 
we cluster all centers $\{\bm_j: j\in L(q+\delta)\}$
such that $\<\bm_j,\bm_{j'}\>/n \ge q$. Each thus formed
cluster $\cC_{\ell}$, $\ell\in L(q)$, is associated
to a new center $\bm_{\ell}$ that is defined as the barycenter of the 
$\{\bm_j:j\in\cC_{\ell}\}$. 

This tree is well defined by ultrametricity, it is roughly balanced, and its 
root will be associated to $\bm_0\approx \bzero$. It is possible that at a certain level $q$
the construction is trivial in the sense that the clusters there consist of a single 
center at level $q+\delta$, in which case level $q$ can be skipped.
Two important properties of this construction are
\begin{align}
\frac{1}{n}\|\bm_{\ell}\|^2 \approx q\;\; \forall \ell\in L(q)\, ,
\;\;\;\;\; \frac{1}{n}\<\bm_{\ell},\bm_{j}-\bm_\ell\>=0\;\; 
\forall j\in \cC_{\ell}\, .
\end{align}
In other words, descending a branch of the tree, one
proceeds from the origin towards the surface of the ball of
radius $\sqrt{n}$,
via a sequence of orthognal increments.

Motivated by this picture, we rename the iterations of
the AMP algorithm   \eqref{eq:iamp-m} as 
$z^{t}$ and $m^{t}$ with $t\in\{0,\delta,\dots, 1-\delta,1\}$,
choose the update functions $f_t$ such that
\begin{align}
  \frac{1}{n}\,\< \bm^{t+\delta}-\bm^t,\bm^t\> &= o_n(1), \label{eq:ortho}\\
  \frac{1}{n}\,\|\bm^{t+\delta}-\bm^t\|^2 &= \delta, \quad \delta\to0 \text{ after } n\to\infty. \label{eq:delta}
\end{align}
A convenient way to enforce the orthogonality constraint is to
choose $f_t$ to be of the form
\begin{align}
  f_t(z^{\le t}) = \sum_{s\in \delta\naturals\cap[0,t-\delta]}
  u(s,z^{\le s})\cdot(z^{s+\delta}-z^{s})\, .
\end{align}
Then state evolution implies the existence of random variables 
$(Z_t,U_t: \, t\in\naturals\delta)$ with the $U_t$ neasurable on 
the $\sigma$-algebra generated by $(Z_s, s\le t)$ and $M_t$
that are the asymptotic equivalent of $\bz^t,  u(t,\bz^{\le })$.
These verify
\begin{align}
&M_t = \sum_{s\in \delta\naturals\cap[0,t-\delta]}\!\!
U_s\cdot(Z_{s+\delta}-Z_{s}),\\
& \E\{Z_{t+\delta}Z_{s+\delta}\}=
\xi'\big(\E\{M_tM_s\}\big)\, ,\label{eq:SE-IAMP}
\end{align}
with $Z_0=0$ and the $(Z_t: t\ge 0)$ Gaussian. 
In particular, increment orthogonality  for the $Z$ process up to time $t$
implies that $Z$ is a Gaussian martingale, and therefore 
$M$ also has orthogonal increments up to time $t$, whence 
by Eq.~\eqref{eq:SE-IAMP}, we get increment orthogonality for $Z$
up to time $t+\delta$.

\subsection{Computational hardness}

The proof of Theorem \ref{thm:alg-guarantee} also shows that 
$\sP_{\cuL}^*(\xi) = \inf_{\gamma\in\cuL} \sP(\gamma;\xi)$ 
is the maximum value of $H(\bsigma)$ achievable by an IAMP algorithm 
(under the assumption that 
  $\inf_{\gamma\in\cuL} \sP(\gamma;\xi)$ is achieved).

Remarkably, $\sP_{\cuL}^*(\xi)$ is also 
the optimal value achievable by a broader class of methods,
known as `overlap concentrated algorithms'
\cite{huang2025tight}. 
In order to define the latter, we will say that 
two Hamiltonians $H^{(1)}(\, \cdot\, ), H^{(p)}(\,\cdot\,)$ 
each with the distribution defined by Eq.~\eqref{eq:MixedPSpinModel}
are $p$-correlated if they are jointly Gaussian with covariance
\begin{align}
\E\{H^{(1)}(\bsigma_1)H^{(p)}(\bsigma_2)\} = 
n\sqrt{p}\, \xi\Big(\frac{\<\bsigma_1,\bsigma_2\>}{n}\Big)\, .
\end{align}
Equivalently, we have $H^{(p)}(\bsigma) =
 \sqrt{p}H^{(1)}(\bsigma) + \sqrt{1-p}\, \tilde{H}^{(1)}(\bsigma)$ for 
 $\tilde{H}^{(1)}(\, \cdot\, )$
an independent copy of $H^{(1)}(\,\cdot\,)$. 
Further, we will identify algorithms with sequences of 
functions $\cA_n$ (indexed by $n$) that take as input a 
Hamiltonian $H:\Sigma_n\to \reals$ and output a point in the convex 
hull of the set  of spin configurations: $\cA_n: H\mapsto
\cA(H)\in \conv(\Sigma_n)$. (We will drop the dependence of $\cA_n$ on $n$.)
\begin{definition}
 We say that an algorithm $\cA: H\mapsto
\cA(H)\in \conv(\Sigma_n)$ is 
\emph{exponentially overlap concentrated} if,
for all $\lambda>0$, there exists $c_0(\lambda)>0$
such that, for all $p\in [0,1]$, we
have (for $H^{(1)}$, $H^{(p)}$ $p$-correlated)
\begin{align}
\prob\Big\{\Big|\<\cA(H^{(1)}),\cA(H^{(p)})\>-
\E\<\cA(H^{(1)}),\cA(H^{(p)})\> \Big|\ge n\lambda\Big\} \le 2\,
e^{-nc_0(\lambda)}\, .\label{eq:overlap-concentrated}
\end{align}
\end{definition}

\begin{remark}
The definition of `overlap concentrated' algorithms
was introduced in \cite{huang2025tight}.  The original definition did 
not require the probability in Eq.~\eqref{eq:overlap-concentrated} 
to be exponentially small. Here we use this
stronger condition since it results in a slightly simpler statement.
\end{remark}

While the above definition only requires 
the algorithm output to be $\cA(H)\in \conv(\Sigma_n)$,
it can be shown \cite{sellke2021approximate} that, 
given $\bm:=\cA(H)\in\conv(\Sigma_n)$,
we can find (for any $\eps>0$) $\bsigma\in\Sigma_n$ `close' to $\bm$
such that --with high probability-- $H(\bsigma)\ge H(\bm)-n\eps$.

Identifying a Hamiltonan $H(\,\cdot\,)$ with the 
vector of Gaussian coefficients $\bg = 
(G^{(k)}_{i_1,\dots, i_k}: 2\le k\le k_{\max},
i_1,\dots,i_k\le n)$, cf. Eq.~\eqref{eq:MixedPSpinModel},
 we can  write $\cA(H)=\cA_G(\bg)$ for 
for some function $\cA_G: \reals^{D}\to \conv(\Sigma_n)$
(with $D= n^2 + \cdots + n^{k_{\max}}$).
Note that $\bg$ is a standard Gaussian vector.
By the Borell-TIS inequality, a sufficient condition for
$\cA$ to be exponentialy overlap concentrated is that
$\bG\mapsto\cA_G(\bg)$ is $L$-Lipschitz 
for some $n$-independent $L$. Several standard optimization algorithms 
can be shown to be Lipschitz.

The following result was proven in \cite{huang2025tight} 
and shows that exponentially overlap concentrated algorithms cannot surpass the 
barrier  $\sP^*_{\cL}(\xi)$.
\begin{theorem}\label{thm:overlap-concentrated}
  Let $\cA: H\mapsto
\cA(H)\in \conv(\Sigma_n)$ be an exponentially overlap
 concentrated algorithm. Then for all $\eps>0$, there exists 
 $c_*(\eps)>0$ such that,
 \begin{align}
  \prob\Big(\frac{1}{n}H(\cA_n(H))\ge \sP^*_{\cL}(\xi)+\eps\Big) \le 2\,
  e^{-nc_*(\lambda)}\, .
 \end{align}
\end{theorem}

This result belongs to a line of work that establishes
computational obstructions by leveraging 
specific geometric properties of the set of near optima
(generically referred to as `overlap gap property') \cite{gamarnik2014limits}.
The geometric property used in the proof of Theorem 
\ref{thm:overlap-concentrated} is motivated by
 the statistical physics picture of the space of near optima.
 
 Additional useful intuition comes from
 the mechanism exploited by the IAMP algorithm 
 of the previous section 
 \cite{subag2021following,montanari2019optimization,alaoui2020optimization}.
We observe that, when one of such algorithms achieves a certain value 
$H(\bsigma^{\salg})\approx n u_*$, then it can be modified to 
generate two configurations $\bsigma^{\salg}_1$ and $\bsigma^{\salg}_2$ 
such that $H(\bsigma^{\salg}_1)\approx H(\bsigma^{\salg}_2)\approx n u_*$
and   $\<\bsigma^{\salg}_1,\bsigma^{\salg}_2\>\approx n t_*$,
for any\footnote{For the algorithm of
\cite{montanari2019optimization,alaoui2020optimization}
 we can run the algorithm as described in the previous section up to
  time $t_*$ (i.e. for $t_*/\delta$
iterations), thus generate $\bm^{t_*}$ with $\|\bm^{t_*}\|^2/n\approx t_*$.
We then set $\bm^{t_*+\delta}_1= \bm^{t_*}+\sqrt{\delta}\bw_{1}$
and $\bm^{t_*+\delta}_2 = \bm^{t_*}+\sqrt{\delta}\bw_{2}$
for $\bw_{1}$ and $\bw_{2}$ two independent standard Gaussian vectors.
Finally, one can apply the IAMP iteration after time $t_*+\delta$.}
$t_*\in (0,1)$.
By a repeated application of the same argument, for any $K\in\naturals$, 
we can generate
configurations $\bsigma^{\salg}_1$, \dots,  $\bsigma^{\salg}_K$ 
such that $H(\bsigma^{\salg}_i)\approx n u_*$ for all $i\le K$,
and the Gram matrix 
$\hbQ_n:=(\<\bsigma^{\salg}_i,\bsigma^{\salg}_j\>/s)_{i,j\le K}$
is $\delta$-close to any matrix with $Q_{ii}=1$ and satisfying
the ultrametric inequality $Q_{i,j}\ge \max(Q_{i,k},Q_{k,j})$ 
for all $i,j,k\le K$ (namely $\|\hbQ_n-\bQ\|_F \le \delta$).
We denote this set of configurations by $\Sigma_{n,k}(\bQ;\delta)
\subseteq (\Sigma_n)^k$.

The bound $H(\bsigma^{\salg})/n\le \sP^*_{\cL}(\xi)+\eps$
then follows from the condition that a $K$-uple of such configurations
$\bsigma^{\salg}_1$, \dots,  $\bsigma^{\salg}_K$ exists, where $\bQ$ 
is the distance matrix of the leaves $L(\TT)$ of regular $b$-ary tree $\TT$
with a large number of levels (in particular $K=b^{m}$ for some 
$m\in\naturals$).

While this discussion provides a geometric interpretation
of the barrier  $\sP^*_{\cL}(\xi)$, and explains why certain algorithms should
fail beyond this point, it does not explain why all exponentially overlap concentrated 
algorithms should fail.
The idea is construct $K$ Hamiltonians assoviated with the leaves
of $\TT$, $\{H_{(i)}: i\in L(\TT)\}$, such that 
$H_{(i)}(\,\cdot\,)$ and  $H_{(j)}(\,\cdot\,)$ are $p_{|i\wedge j|}$
correlated. Here $p_0\le p_1\le , \dots, \le p_m$ are $m$ parameters and
$|i\wedge j|$ is the level of the least common ancestor of $i$ and $j$.
By properly choosing the parameters $p_k$ (depending on the algorithm),
we can ensure that the outputs $\bsigma_i^{\salg}$ of the algorithm on inputs 
$H_{(i)}$, $i\in L(\TT)$ is such that the resulting Gram matrix 
approximates $\bQ$. 

Finally, the proof is completed by constructing $\TT$ and $\bQ$ 
such that (for $\delta$ sufficiently small), with high probability,
\begin{align}
\frac{1}{n|L(\TT)|}\max_{(\bsigma_i) \in \Sigma_{n,k}(\bQ;\delta)} 
\sum_{i\in L(\TT)}
H_{(i)}(\bsigma_i) \le \sP^*_{\cL}(\xi)+\eps\, .
\end{align}
In particular, the upper bound is proven by Guerra's
interpolation method \cite{guerra2003broken,aizenman2003extended}
(see \cite{panchenko2007overlap} for MT's work on a related question).

\section{A broader survey}
\label{sec:Broader}

Techniques developed for the study of the quadratic optimization problem 
\eqref{eq:Hsigma} (for $\bA\sim\GOE(n)$) or its the mixed $p$-spin model 
\eqref{eq:MixedPSpinModel} have been generalized to a broad array of other
 problems.
We provide a few examples and pointers below, while keeping in mind that 
our mathematical understanding is not equally advanced in each of these cases.

\paragraph{Tensor models, a.k.a. multi-spin models.}
As mentioned already, the mixed $p$-spin model \eqref{eq:MixedPSpinModel}
has attracted a significant amount of work within statistical physics.
An important motivation to move beyond the quadratic
cost \eqref{eq:Hsigma} lies in the fact that the mixed $p$-spin model 
can present (for certain choices of the coefficients $\xi_k$)
new phase transitions that are not observed in the original SK model.
These phase transitions have been used as toy models for 
phase transitions in structural glasses 
\cite{kirkpatrick1987p,crisanti1992spherical}.

Within computer science, the optimization problem $\max_{\bsigma\in\Sigma_n}H(\bsigma)$ has been mainly studied in
the case $\Sigma_n = \sqrt{n}\S^{n-1}$, for the homogeneous case 
$H(\bsigma) =\<\bG^{(k)},\bsigma^{\otimes k}\>/n^{(k-1)/2}$. In this case,
it is equivalent to computing the injective norm of the tensor $\bG^{(k)}$,
a notoriously hard problem for 
$k\ge 3$ \cite{barak2012hypercontractivity,bhattiprolu2017sum}.

A number of statistical problems have been studied in recent years 
for which the data take the form of a high dimensional tensor
 $(\bY\in\reals^n)^{\otimes k}$.
The simplest example (introduced in \cite{montanari2014statistical}) is 
`tensor PCA': we observe a tensor 
$\bY = \lambda\btheta^{\otimes k}+\bG^{(k)}/\sqrt{n}$, with $\bG$ as above, $\btheta$
an unknown unit vector, and $\lambda$ a fixed signal-to-noise ratio. 
We would like to estimate $\btheta$ from data $\bY$. 

\paragraph{Generalized perceptrons.} Given $h:\reals\times\reals\to \reals$, and random variables
$\bx_1,\dots, \bx_n\in\reals^d$, $y_1,\dots, y_n\in\reals$,
consider  (for $\bsigma\in \Sigma_d\subseteq\reals^d$)
\begin{align}
H(\bsigma) =\sum_{i=1}^n h(\<\bx_i,\bsigma\>, y_i)\, .\label{eq:GenPerceptron}
\end{align}
A well studied model for the randomness assumes $\bx_i\sim\normal(0,\id_d/d)$.

The most studied example in spin glass theory assumes  
$h(z,y) = h_0(z)$ independent of $y$, with
$h_0(z) < 0$ for $z< \kappa$ and $h_0(z) = 0$ for $z\ge \kappa$,
and either $\Sigma_d= \sqrt{d}\S^{d-1}$ (spherical perceptron \cite{gardner1988space}) 
or  $\Sigma_d= \{-1,1\}^d$  (Ising perceptron \cite{krauth1989storage}). 
In this case a maximizer $\bsigma$ with $H(\bsigma)=0$
corresponds to a vector $\bsigma\in \Sigma_d$ in the intersection of 
half-spaces $\cap_{i\le n}\{\bsigma\in\reals^d:\, \<\bsigma,\bx_i\>\ge \kappa\}$.

In one of his earliest rigorous contributions to spin glass theory,
MT proved that the cardinality of set of minimizers  $\cuZ_{n,d} := 
\Sigma_d\cap
_{i\le n}\{\bsigma\in\reals^d:\, \<\bsigma,\bx_i\>\ge \kappa\}$
concentrates on the logarithmic scale. Namely,
$n^{-1}\log|\cuZ_{n,d}|= \phi_{n,d}+o_P(1)$ for some deterministic quantity 
$\phi_{n,d}=O(1)$
\cite{talagrand1999self}.

This problem has attracted renewed interest in theoretical computer science
because of its connection to discrepancy theory.
Given a matrix $\bX\in\reals^{n\times d}$, a typical question in discrepancy theory
is to find the minimum value of  $\|\bX\bsigma\|_p$ over $\bsigma\in\{+1,-1\}^d$.
\cite{spencer1985six,bansal2010constructive}. For $\bX$ a random matrix,
 this problem fits the above formulation
with the choice $h(z,y)= -|z|^p$.

Finally, the objective function \eqref{eq:GenPerceptron} plays an important 
role in high-dimensional statistics and machine learning. Indeed,
understanding the `space of interactions of neural network models'
was Elizabeth Gardner's original motivation for introducing this problem
 in \cite{gardner1988space}.

 In a canonical model 
for supervised learning, the pairs 
$(\bx_i,y_i)$, $i\le n$ are  assumed to be i.i.d. from some unknown distribution
on $\reals^d\times \reals$. The vectors, $\bx_i\in\reals^d$ are referred to as `covariates' or 
`feature vectors',
while the $y_i\in \reals$ are `responses' or `labels'. To give a concrete example $\bx_i$ 
could encode the features of a house on sale, and $y_i$ the sales price.
We would like to learn a model $f:\reals^d\to\reals$ to predict the
response $y_{\snew}$ associated to a new feature vector $\bx_{\snew}$.

A popular approach to achieve this is to minimize the error on the
`training data':
\begin{align}
\mbox{minimize}&\;\;\;\; \sum_{i=1}^n \big(y_i-f(\bx_i)\big)^2\ ,\\
\mbox{subj. to}&\;\;\;\;  f\in \cuF\, ,
\end{align}
for $\cuF$ a suitable function classes (of course the mean square error can be replaced by other 
metrics). If we consider `generalized linear models', of the form 
$\cuF_{\varphi}:= \{f(\bx) = \varphi(\<\bsigma,\bx\>):\; \|\bsigma\|=\sqrt{d}\}$,
for a fixed `activation function' $\varphi$,
then this cost function reduces to a special case of of the Hamiltonian
$H(\bsigma)$ of \eqref{eq:GenPerceptron}.

We refer to  \cite{barbier2019optimal,ding2025capacity,nakajima2023sharp,montanari2024exceptional} for a few recent pointers to the 
literature on this model, and its generalizations.

\paragraph{Random $k$ satisfiability.} Informally, we focused so far on models with `dense' interactions.
In the SK model, the Hamiltonian takes the form \eqref{eq:Hsigma}, whereby 
$\bA\sim\GOE(n)$, hence any two spins $\sigma_i$, $\sigma_j$ interact via the term 
$A_{ij}\sigma_i\sigma_j$ which is of order $1/\sqrt{n}$. Similarly, in the multi-spin
model \eqref{eq:MixedPSpinModel}, every $k$-uple of spins interacts via 
$G^{(k)}_{i_1,\dots,i_k}\sigma_{i_1}\cdots \sigma_{i_k}/n^{(k-1)/2}$, which is of order
$n^{-(k-1)/2}$. Finally, in the generalized perceptron of Eq.~\eqref{eq:GenPerceptron},
spin $\sigma_j$ contribution to the interaction $h(\<\bx_i,\bsigma\>,y_i)$
is $x_{ij}\sigma_j$, which of the same order $n^{-1/2}$.

The fact that all interactions tems are of the same order
implies certain simplifications. To understand why, we can compute
 the energy change due to flipping spin $\sigma_i$,
 $\Delta_i(\bsigma):=H(\bsigma^{(i)})-H(\bsigma)$
 (here $(\bsigma^{(i)})_i=-\sigma_i$,
while $(\bsigma^{(i)})_j=\sigma_j$ for all $j\neq i$.)
 For simplicity, consider the case in which
$\bsigma\in\{+1,-1\}^n$ is independent of the Hamiltonian $H(\,\cdot\,)$. 
 It is not hard to see that 
$\Delta_i(\bsigma)$  is approximately Gaussian 
under  any of the above models for $H(\,\cdot\,)$ (it is 
sufficient that the coefficients $G^{(k)}_{i_1\dots i_k}$ are i.i.d. 
centerd, under finite second moment, but not necessarily Gaussian).
Fr fense models,
such Gaussian approximations are possible for suitably defined quantities
even if $\bsigma$ is configuration from the Gibbs measure,

Starting in the 1980s, spin glass theory techniques have been extended 
to models
with `sparse' interactions in which such simplifications do not apply.
The earliest example is the Viana-Bray model \cite{viana1985phase} which 
again takes the form
\eqref{eq:Hsigma}, with $\bA=\bA^{\sT}$ a sparse random matrix,
with $A_{ij}\neq 0$ only on if $(i,j)\in E$, with $E$ the edge set of an 
Erd\"os-Renyi random graph, with average degree of order one. 
For instance, we can
take $\prob(A_{ij}=+1) = \prob(A_{ij}=-1)  =\gamma/2n$ and  $\prob(A_{ij}=0) = 1-\gamma/n$
for some constant $\gamma$.
(See  \cite{mezard2001bethe,panchenko2004bounds,el2023local} for further 
results on this and related models.)

The generalization of spin glass techniques 
to models on sparse graph
has attracted significant interest because of its role in
 the analysis of random combinatorial 
optimization problems. A prototypical
example is random $k$-satisfiability ($k$-SAT). 
An instance of the $k$-SAT problem requires to assign values to 
$n$ Boolean variables $\sigma_1,\dots,\sigma_n$ as to satisfy 
$m$ constraints. Each constraint (`clause') is the logical OR of $k$ 
variables or their negations.
Equivalently, the variables assignment must satisfy a Boolean formula which
is the logical AND of  $m$ clauses. An example of such a formula
with $k=3$, $n=5$, $m=4$ is (here overline denotes logical negation)
\begin{align}
(\sigma_1\vee \osigma_3\vee \sigma_5)\wedge (\sigma_2\vee \sigma_3\vee \osigma_4)
\wedge (\sigma_2\vee \osigma_4\vee \osigma_5)\wedge (\sigma_1\vee \sigma_3\vee \sigma_4)\, .
\end{align}
Of course we can encoded Boolean variables as $ \sigma_i\in\{+1,-1\}$ and 
denote by $E_a(\bsigma)$ the indicator function for clause $a$ being violated
 by $\bsigma\in\{+1,-1\}^n$. Finding a satisfying assignment thus amounts 
 to minimizing the Hamiltonian
\begin{align}
H(\bsigma) = \sum_{a=1}^m E_a(\bsigma)\, .
\end{align}
This Hamiltonian is always non-negative
and formula is satisfiable if and only if $\min_{\bsigma\in\Sigma_n} 
H(\bsigma)=0$ or,
in physics jargon, the ground state energy is equal to zero. 
The problem of deciding whether a formula is
satisfiable or not (i.e. whether the ground state energy is zero or strictly 
positive) is referred to as a `decision problem' 
in theoretical computer science\footnote{From a computational complexity perspective,
decision and optimization are roughly equivalent.}.

Since  $k$-SAT is intractable in the worst case (NP-complete), 
computer scientists have
introduced probability distributions over random formulas as a way 
to compare different heuristics \cite{gu1997algorithms}.
`Random $k$-SAT' is the name of one particularly simple (and useful)
 such distributions, whereby 
a formula is drawn uniformly at random among $k$-SAT all formulas with $m$ clauses and $n$
variables (equivalently, each of the $m$ clauses is drawn uniformly at random among all
the $2^k\binom{n}{k}$ possible $k$-clauses). Of particular interest is limit of
large formulas $m,n\to\infty$ with $m/n\to \alpha\in (0,\infty)$.

Physics approaches study a `finite-temperature' version of this cost function whereby 
 one introduces the partition function 
 $Z_n(\beta) = \sum_{\bsigma\in\Sigma_n}e^{-\beta H(\bsigma)}$.
 The limit $\beta\to\infty$ of the free energy density allows to characterize
 the asymptotics of the minimum fraction of unsatisfied clauses:
 \begin{align*}
 e_*(\alpha) &:= \lim_{n\to\infty} 
 \frac{1}{n}\min_{\bsigma\in\{+1,-1 \}^n} H(\bsigma) 
 =\lim_{\beta\to\infty}\lim_{n\to\infty} \frac{1}{n\beta}\log Z_n(\beta)\,.
 \end{align*}
 Limits above were shown to exist in expectation (and hence almost surely
 by simple concentration arguments) in  \cite{franz2003replica}
 by extending Guerra-Toninelli's interpolation argument 
 \cite{guerra2002thermodynamic}. 
In collaboration with Panchenko (and concurrently with Franz, Leone, Toninelli 
\cite{franz2003replica,franz2003replica_b}), MT 
adapted Guerra's interpolation approach to this problem establishing that 
physicists
predictions for $\lim_{n\to\infty} n^{-1}\log Z_n(\beta)$ provided indeed 
rigorous upper bounds
\cite{panchenko2004bounds}.

Much remains to be understood in the study of models on spares graphs, especially
on algorithmic aspects. We refer to
 \cite{monasson1999determining,mezard2001bethe,mezard2002analytic,
 ding2016maximum,coja2017information,ding2022proof,el2023local} 
 for a few recent pointers to the literature on sparse graph models.

\section*{Acknowledgements}

AM was partially supported by the NSF through award DMS-2031883, the Simons Foundation through
Award 814639 for the Collaboration on the Theoretical Foundations of 
Deep Learning.

\newcommand{\etalchar}[1]{$^{#1}$}
\providecommand{\bysame}{\leavevmode\hbox to3em{\hrulefill}\thinspace}
\providecommand{\MR}{\relax\ifhmode\unskip\space\fi MR }
\providecommand{\MRhref}[2]{%
  \href{http://www.ams.org/mathscinet-getitem?mr=#1}{#2}
}
\providecommand{\href}[2]{#2}

\addcontentsline{toc}{section}{References}

\end{document}